\documentclass[11pt]{amsart}

\usepackage{enumerate}
\usepackage{latexsym}
\usepackage[pdftex]{graphicx}
\usepackage{amssymb}
\usepackage{xcolor}
\usepackage{stackrel}
\usepackage[normalem]{ulem}

\newtheorem{theorem}{Theorem}
\newtheorem{lemma}{Lemma}
\newtheorem{proposition}{Proposition}
\newtheorem{remark}{Remark}
\newtheorem{example}{Example}

\newtheorem{definition}{Definition}
\newtheorem{corollary}{Corollary}

\newtheorem{problem}{Problem}
\newtheorem*{problem*}{Problem}

\newcommand{\beq}{\begin{equation}}
\newcommand{\eeq}{\end{equation}}
\newcommand{\beqna}{\begin{eqnarray*}}
\newcommand{\eeqna}{\end{eqnarray*}}
\newcommand{\beqn}{\begin{equation*}}
\newcommand{\eeqn}{\end{equation*}}
\newcommand{\bp}{\begin{proof}}
\newcommand{\ep}{\end{proof}}
\newcommand{\bprop}{\begin{proposition}}
\newcommand{\eprop}{\end{proposition}}
\newcommand{\bt}{\begin{theorem}}
\newcommand{\et}{\end{theorem}}
\newcommand{\bex}{\begin{example}}
\newcommand{\eex}{\end{example}}
\newcommand{\bc}{\begin{corollary}}
\newcommand{\ec}{\end{corollary}}
\newcommand{\bl}{\begin{lemma}}
\newcommand{\el}{\end{lemma}}
\newcommand{\bprob}{\begin{problem}}
\newcommand{\eprob}{\end{problem}}
\newcommand{\br}{\begin{remark}}
\newcommand{\er}{\end{remark}}
\newcommand{\bd}{\begin{definition}}
\newcommand{\ed}{\end{definition}}

\begin{document}

\title
[On an equichordal property of a pair of   convex bodies]
{On an equichordal  property of a pair of  convex bodies}

\author[D. Ryabogin]{Dmitry Ryabogin}
\address{Department of Mathematics, Kent State University,
Kent, OH 44242, USA} \email{ryabogin@math.kent.edu}

\thanks{The   author is  supported in
part by  Simons Collaboration Grant for Mathematicians program 638576 and by U.S.~National Science Foundation Grant
DMS-1600753.}

\keywords{convex  bodies, equichordal property}

\begin{abstract}
Let $d\ge 2$ and let $K$ and $L$ be two  convex bodies in ${\mathbb R^d}$ such that $L\subset \textrm{int}\,K$ and the  boundary of $L$ does not contain a segment.  If
$K$ and $L$ satisfy  the $(d+1)$-equichordal property, i.e.,   
   for any line $l$ supporting  the boundary of $L$ and the points $\{\zeta_{\pm}\}$  of the intersection of  the boundary of $K$ with $l$,
$$
\textrm{dist}^{d+1}(L\cap l, \zeta_+)+\textrm{dist}^{d+1}(L\cap l, \zeta_-)=2\sigma^{d+1}
$$
holds, where    the  constant $\sigma$ is independent of $l$, does it follow that $K$ and $L$ are concentric Euclidean balls?
We   prove that 
 if  $K$ and $L$ have $C^2$-smooth boundaries and 
  $L$ is a  body of revolution,
 then $K$ and $L$ are 
 concentric Euclidean balls. 
\end{abstract}

\maketitle

\section{Introduction}

\begin{figure}[ht]
	\centering
	\includegraphics[height=3in]{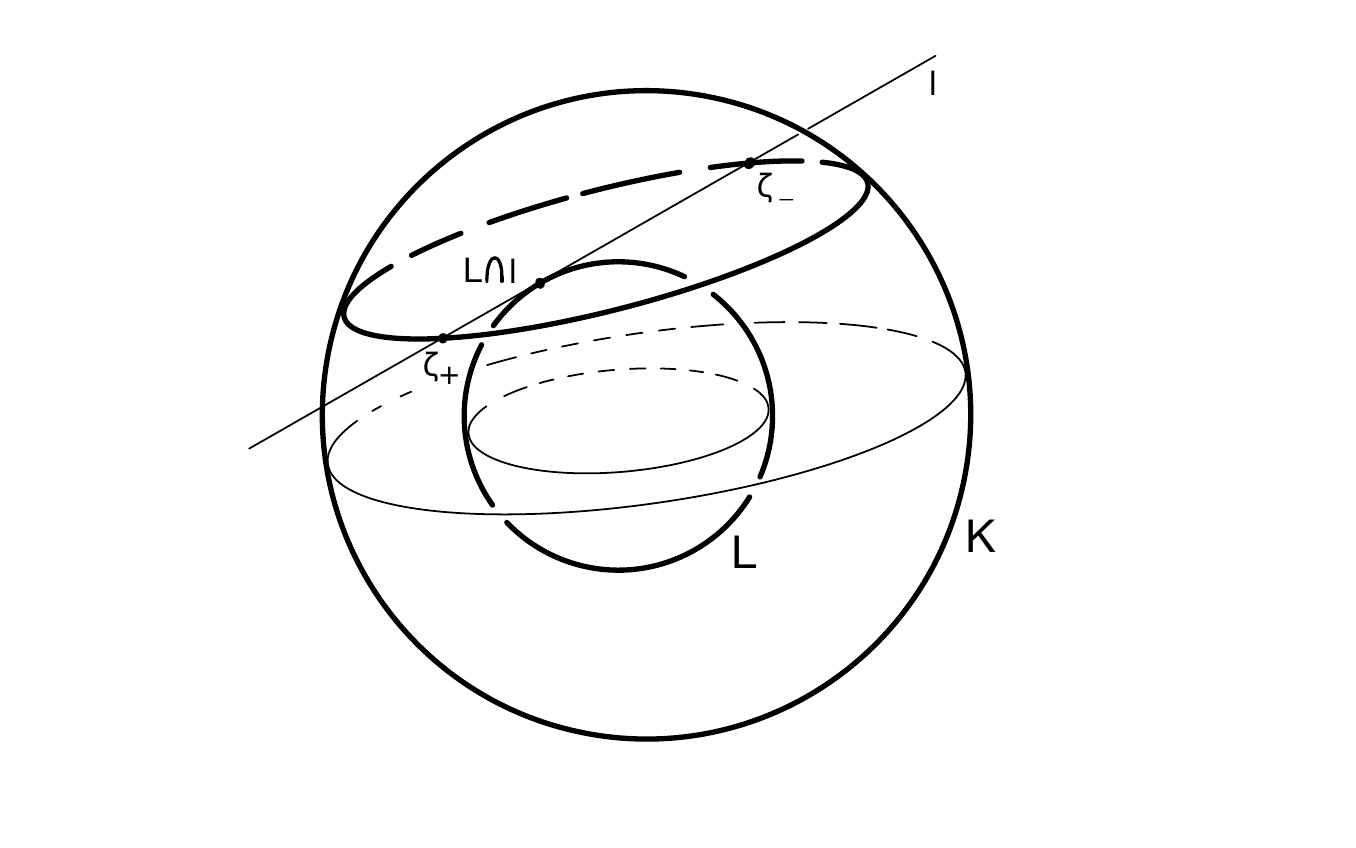} 
	\caption{We have $\textrm{dist}^i(L\cap l, \zeta_+)+\textrm{dist}^i(L\cap l, \zeta_-)=2\sigma^i$.}
	\label{fonseca}
\end{figure}

Let $d\ge 2$
and let $K$ and $L$ be two convex bodies in ${\mathbb R^d}$
such that $L\subset \textrm{int}\,K$ and the  boundary of $L$ does not contain a segment. For any line $l$ supporting $L$ we consider two points $\zeta_{\pm}$ of the  intersection of the boundary of $K$ with $l$.
Given $i\in {\mathbb R}$ we say that the bodies $K$ and $L$ satisfy  the $i$-equichordal property if   there exists a constant $\sigma$  independent of $l$ such that
\begin{equation}\label{hop1}
\textrm{dist}^i(L\cap l, \zeta_+)+\textrm{dist}^i(L\cap l, \zeta_-)=2\sigma^i,
\end{equation}
(see Figure \ref{fonseca}).
If $i=0$ we replace (\ref{hop1}) with
\begin{equation}\label{hop2}
\textrm{dist}(L\cap l, \zeta_+)\textrm{dist}(L\cap l, \zeta_-)=\sigma^2,
\end{equation}
(cf. \cite{Ga}, page 233).

\bprob\label{s1}
 Let $d\ge 2$ and $i\in{\mathbb R}$. Are two concentric Euclidean balls the only pair of bodies in ${\mathbb R^d}$  satisfying  the $i$-equichordal property? 
\eprob

Similar questions to that of the problem above   were  raised in \cite{Sa},  
\cite{BaL},  \cite[A1, page 9]{CFG},  \cite{YZ}; see also \cite{RYZ} and references therein.
In particular, it is known that the answer to Problem \ref{s1} is affirmative for $d\ge 3$, provided
 $L$ is a Euclidean  ball, \cite{BaL}.

We would also like to  mention several results related to the  connection between Problem \ref{s1}  and 
Problem 19 of Ulam from the Scottish book, which asks
{\it  if a solid of uniform density which floats in water in every position is necessarily a sphere}, \cite[page 90]{M},  \cite[A9, page 19]{CFG}. 

The  plane counterexamples to Ulam's problem constructed in \cite{A}, \cite{Weg1}, \cite{Weg2},
show that for $d=2$, $i=1$,
the answer  to Problem \ref{s1}  is negative,  even under the additional assumption that for every line $l$ supporting $L$, the point of tangency $L\cap l$ divides the chord $K\cap l$  into two parts of equal length. On the other hand, it is known \cite{BMO} that, under this division assumption and under the assumption that $l$ divides the boundary of $K$ in constant ratio $\frac{\mu(\sigma)}{1-\mu(\sigma)}$ for $\mu=\frac{1}{3}$, 
 $\mu=\frac{1}{4}$, $\mu=\frac{1}{5}$, and $\mu=\frac{2}{5}$,
 the answer to Problem \ref{s1} is affirmative; see also \cite{Od}. 
Additionally,  if $d\ge 3$,  and if  for
 every line $l$ supporting $L$ the point 
$L\cap l$ divides the chord $K\cap l$ into two parts of equal length, then the answer to Problem \ref{s1} is affirmative,  \cite{O}. 
Finally, we remark that
a negative answer to Problem \ref{s1} in the case $i=d+1$, $d\ge 3$, 
presents a possibility for a  negative answer to Ulam's conjecture, \cite{R1}, \cite{R2}.

In this paper we prove the following result.

\bt\label{mdauzh}
 Let $d\ge 3$ and let $K$ and $L$ be two convex bodies in ${\mathbb R^d}$ of class  $C^2$  
 satisfying the $(d+1)$-equichordal property. 
 If $L$ is a  body of revolution, then $K$ and $L$ are 
 concentric Euclidean balls.
\et

A similar result  can be proved for  general $i$-equichordal property,  $i\in{\mathbb R}$. Since 
 our interest in Problem \ref{s1} comes, partly,  from its relation to the Problem of Ulam, and since, in our opinion, the proof for $i\ne d+1$ does not add to the ideas when $L$ is a body of revolution,
 we restrict ourselves to the case $i=d+1$.

{\bf Notation and basic definitions. }
 Let $d\ge 2$. 
A convex body $K\subset {\mathbb R^d}$ is a convex compact set with a non-empty interior $\textrm{int}\,K$. We denote by $B^d(r)$ the Euclidean ball centered at the origin of radius $r>0$.
Given  $\xi\in S^{d-1}$ we put
$\xi^{\perp}=\{p\in{\mathbb R^d}:\, p \cdot \xi=0\}$ to be   the subspace orthogonal to $\xi$, and  $p \cdot\xi=p_1\xi_1+\dots +p_d\xi_d$ is the usual inner product in ${\mathbb R^d}$.  

We say that a line $l$ is a supporting line of  a convex body $L$ if $L\cap l\neq \emptyset$, but $\textrm{int}\,L\cap l=\emptyset$.

Let  $m\in{\mathbb N}$. We say that  a convex body $K$ in ${\mathbb R^d}$ is of class  $C^m$ if for every point $z$ on the boundary $\partial K$ of $K\subset {\mathbb R^d}$ there exists a neighborhood $U_z$ of $z$  in ${\mathbb R^d}$ such that  $\partial K\cap U_z$ can be written as a graph of a function having 
all continuous partial derivatives  up to  the $m$-th order.

\section{Auxiliary statements,    $K$ and $L$ are the  bodies of revolution about the same axis in ${\mathbb R^3}$}\label{uxnem}

At first we introduce some convenient notation that helps  to work with bodies of revolution.

\begin{figure}[ht]
	\centering
	\includegraphics[height=3in]{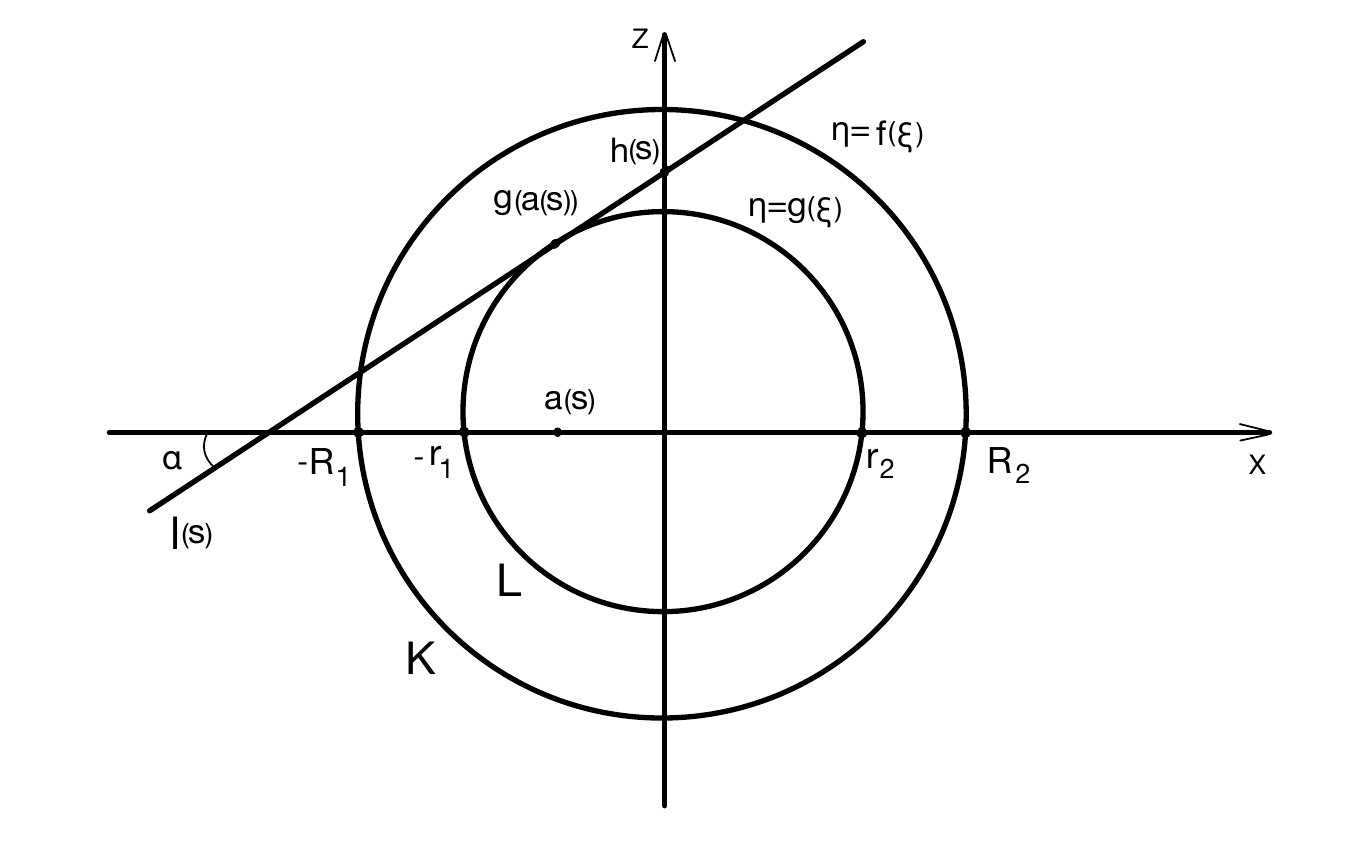} 
	\caption{$K\cap \{(x,y,z):\,y=0\}$ and $L\cap \{(x,y,z):\,y=0\}$ with their boundaries described by the graphs of functions $\eta=f(\xi)$ and $\eta=g(\xi)$.}
	\label{rR}
\end{figure}

Let $K\subset {\mathbb R^3}$ be a body of revolution about the $x$-axis with $C^3$  boundary   described by a function $\eta=f(\xi)\ge 0$ supported by  the segment $[-R_1,R_2]$. Assume also that  $L$ is a body of revolution about the same axis, and its boundary  is described by the function $\eta=g(\xi)\ge 0$ supported by the segment $[-r_1,r_2]\subset (-R_1,R_2)$ (See Figure \ref{rR}).

We will denote by $H_s$  the plane parallel to the $y$-axis and containing the line $l(s)=\{(\xi,0, s\xi+h(s)):\,\xi\in{\mathbb R}\}$,  where $l(s)$ is tangent to the graph of $g$ at the corresponding point $(a(s), 0, g(a(s)))$, $s=\tan\alpha$ with $\alpha\in(-\frac{\pi}{2}, \frac{\pi}{2})$ being the angle between the $x$-axis and $l(s)$, and $h(s)$ is the $z$-intercept of $l(s)$.

Let $s\in {\mathbb R}$ be fixed and let  $\ell_s$ be the line parallel to the $y$-axis passing through $(a(s), 0, g(a(s)))$. Since the section $K\cap H_s$ is symmetric  with respect to the line $l(s)$,  the chord $G_s=K\cap \ell_s$  is divided by $(a(s), 0, g(a(s)))$ into two parts of equal length $\sigma$.

\begin{figure}[ht]
	\centering
	\includegraphics[height=3.5in]{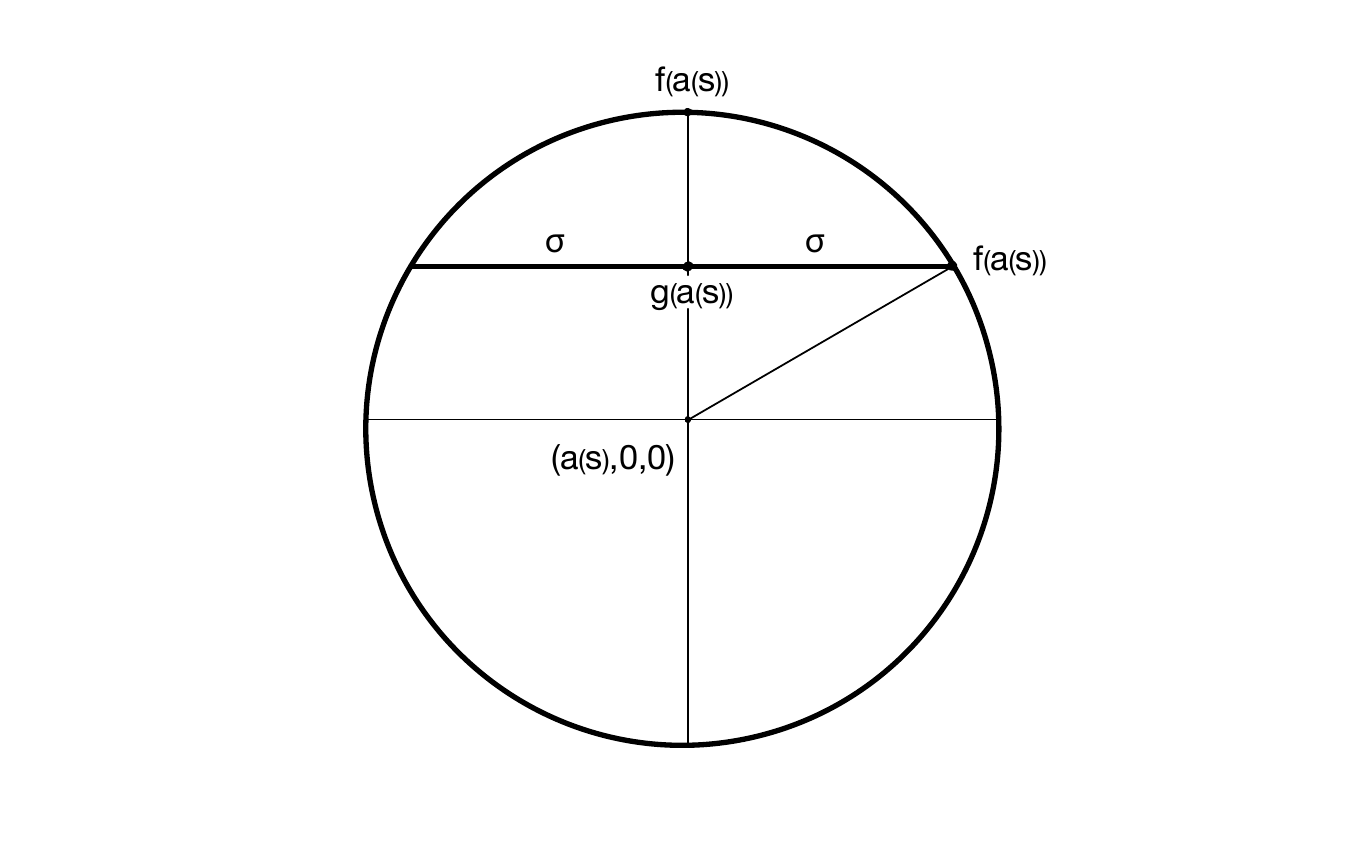} 
	\caption{The section $K\cap\{(x,y,z):\,x=a(s)\}$. We have $g(a(s))=\sqrt{f^2(a(s))-\sigma^2}$.}
	\label{seca}
\end{figure}

Let $s\in {\mathbb R}$ and  $a(s)\in (-r_1,r_2)$. 
 Since  $K$ is a body of revolution, the chord of length $2\sigma$ can be inscribed into a circle of radius $f(a(s))$ only provided
$$
g(a(s))=\sqrt{f^2(a(s))-\sigma^2}
$$
(see Figure \ref{seca}).
Since $K\cap\{ (x,y,z)\in{\mathbb R^3}:\,x=-r_1,r_2\}$ are discs of radius $\sigma$, we have
\begin{equation}\label{tilde}
f(-r_1)=f(r_2)=\sigma, \quad g(\xi)=\sqrt{f^2(\xi)-\sigma^2}\quad \forall \xi\in [-r_1, r_2].
\end{equation}

Translating the bodies if necessary, we can and do assume that $a(0)=0$.

\subsection{Some  results on the $(d+1)$-equichordal plane bodies symmetric with respect to the axis}

Let $P\subset {\mathbb R^2}$ be a convex body  containing the origin in its interior with $C^3$-smooth boundary.  Assume also that $P$ is symmetric with respect to the $x$-axis and  it   satisfies the $(d+1)$-equichordal property with respect to the origin, i.e.,  there exists a constant $\sigma$ such that
$$
\forall \theta\in S^1\qquad \rho_P^{d+1}(\theta)+\rho_P^{d+1}(-\theta)=2\sigma^{d+1},\quad \textrm{for some}\quad d\ge 3.
$$

\begin{figure}[ht]
	\centering
	\includegraphics[height=3in]{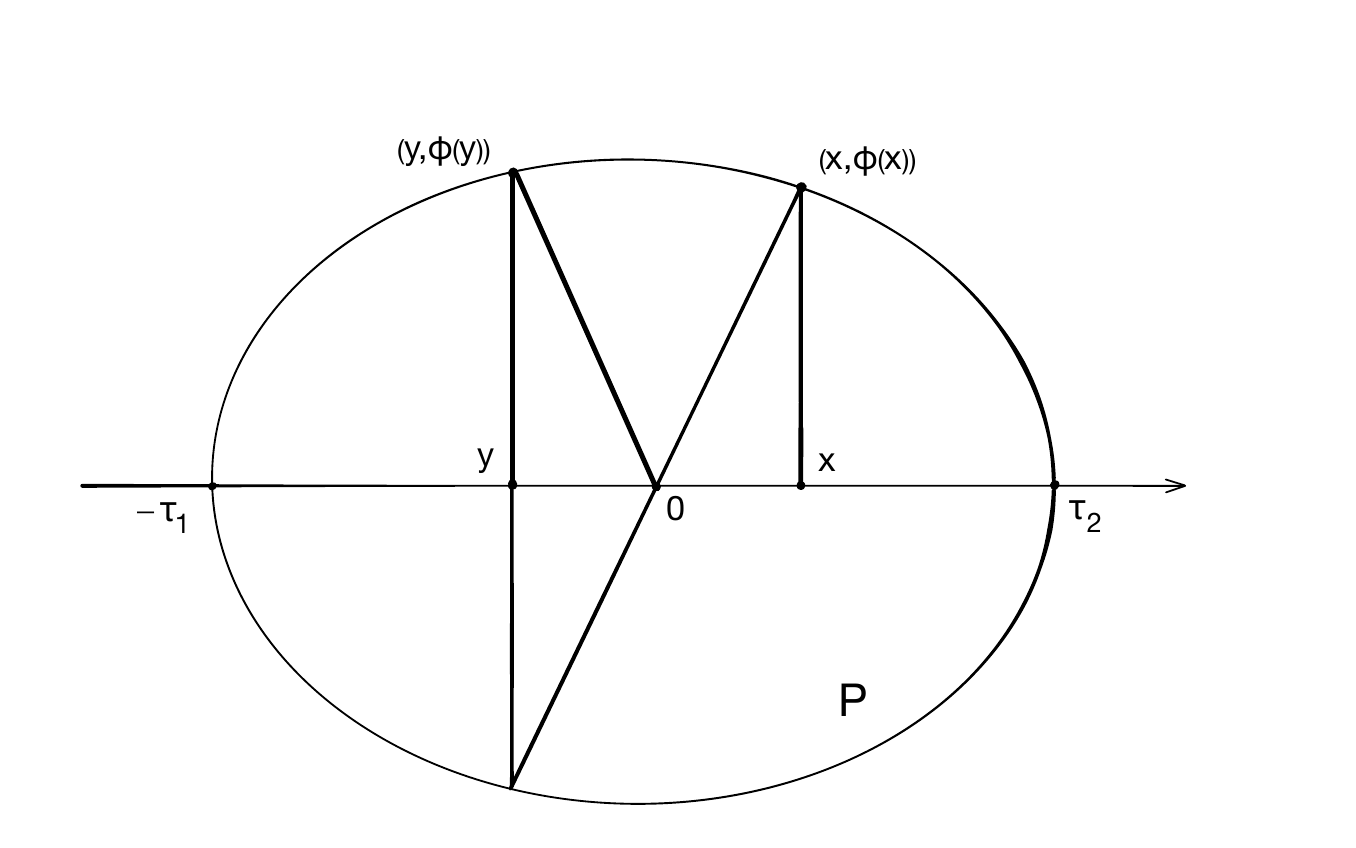} 
	\caption{The equichordiality of $P$.}
	\label{Pu}
\end{figure}

If the upper part of the boundary of $P$ 
is  described by a graph of a positive function $\phi$ on $[-\tau_1,\tau_2]$, then by the Pythagorean Theorem  and the symmetry with respect to the $x$-axis, the function $\phi $  satisfies
\begin{equation}\label{votto}
(x^2+\phi^2(x))^{\frac{d+1}{2}}+(y^2+\phi^2(y))^{\frac{d+1}{2}}=2\sigma^{d+1},
\end{equation}
(see Figure \ref{Pu}).
Here $y\in [-\tau_1,0]$, $\tau_1>0$, and $x\in [0,\tau_2]$, $\tau_2>0$,  are such that
$$
\frac{\phi(x)}{x}=\frac{\phi(y)}{|y|},
$$
i.e., 
\begin{equation}\label{yx}
|y|^{d+1}=\frac{(2\sigma^{d+1}-(x^2+\phi^2(x))^{\frac{d+1}{2}})x^{d+1}}{(x^2+\phi^2(x))^{\frac{d+1}{2}}}.
\end{equation}

To simplify the computations we will write
\begin{equation}\label{chi}
\phi^2(x)=\sigma^2-x^2+\chi(x),\quad x\in [-\tau_1,\tau_2],
\end{equation}
where $\chi$ is a function we want to determine.  By the symmetry of $P$ with respect to the $x$-axis,  we have $\phi(0)=\sigma$, hence, $\chi(0)=0$.  
We rewrite  (\ref{votto}) as
\begin{equation}\label{chu}
(\sigma^2+\chi(x))^{\frac{d+1}{2}}+(\sigma^2+\chi(y))^{\frac{d+1}{2}}=2\sigma^{d+1},\quad y\in [-\tau_1,0],\,x\in [0,\tau_2].
\end{equation}

Let
\begin{equation}\label{per11}
q(z)=\frac{\chi(z)}{\sigma^2}=\frac{\phi^2(z)+z^2-\sigma^2}{\sigma^2}=\frac{\phi^2(z)+z^2}{\sigma^2}-1.
\end{equation}
Then 
conditions (\ref{chu}) and (\ref{yx}) can be  written as
\begin{equation}\label{per}
(1+q(x))^{\frac{d+1}{2}}+(1+q(y))^{\frac{d+1}{2}}=2.
\end{equation}
\begin{equation}\label{yx2}
|y|^{d+1}=\frac{(2\sigma^{d+1}-(\sigma^2+\chi(x))^{\frac{d+1}{2}})x^{d+1}}{(\sigma^2+\chi(x))^{\frac{d+1}{2}}}=\frac{(2-(1+q(x))^{\frac{d+1}{2}})x^{d+1}}{(1+q(x))^{\frac{d+1}{2}}}.
\end{equation}

Our first lemma is technical, but it is crucial for our further considerations.

\bl\label{comp0}
Let $\chi$ be as in (\ref{chu}), where $y$ is as in (\ref{yx2}). Then
\begin{equation}\label{concb}
2 \sigma^2\chi''(0)+(d+1)(\chi'(0))^2=0.
\end{equation}
\el
\bp
Since $q(0)=0$, we can assume that there exists a neighborhood $U_0$ of the origin 
such that $|q(x)|<1$ for all $x\in U_0$.
Using (\ref{yx2}) we see that for $y<0$ and $x>0$ we have
\begin{equation}\label{root1}
|y|=x \,\frac{(2-(1+q(x))^{\frac{d+1}{2}})^{\frac{1}{d+1}}}{(1+q(x))^{\frac{1}{2}}}\qquad \forall x\in U_0.
\end{equation}
We will show at first that 
\begin{equation}\label{Tay1}
|y|=x\Big(1- \varepsilon_1 x+\Big(-\varepsilon_2+\frac{3-d}{4}\varepsilon_1^2\Big) x^2+o(x^2)\Big),\quad \forall x\in V_o,
\end{equation}
where  $V_0\subset U_0$ is a neighborhood of the origin that will be chosen later and $\varepsilon_j$, $j=1,2$, are the Taylor  coefficients  of the decomposition  of $q$ near the origin,
\begin{equation}\label{ep}
q(x)=\varepsilon_1 x+\varepsilon_2 x^2+o(x^2),
\,\,
\varepsilon_j=\frac{q^{(j)}(0)}{j!}=\frac{\chi^{(j)}(0)}{\sigma^2j!},\quad x\in V_0.
\end{equation}

To prove (\ref{Tay1}), we compute the first and  second derivatives of the function
$$
\frac{(2-(1+z)^{\frac{d+1}{2}} )^{\frac{1}{d+1}} }     {(1+z)^{\frac{1}{2}}  }.
$$
Routine calculations show that they are equal to $-1$ and $\frac{3-d}{4}$ respectively, and we can express $y$ via $x$ up to the terms of the second order,
$$
|y|=1-q(x)+\frac{3-d}{4}q^2(x)+o(q^2).
$$
Now we will use (\ref{ep}) and the previous decompositions to obtain
$$
|y|=x\Big(1-(\varepsilon_1 x+\varepsilon_2 x^2)+\frac{3-d}{4}(\varepsilon_1 x+\varepsilon_2 x^2)^2+o(x^2)\Big).
$$
This gives (\ref{Tay1}). 

Next,  we  use  (\ref{Tay1}) to obtain two relations that will lead to 
(\ref{concb}).
We see that 
\begin{equation}\label{sum}
x-|y|=\varepsilon_1 x^2+o(x^2).
\end{equation}
Since 
$$
y^2=x^2(1- 2\varepsilon_1 x)+o(x^3),
$$
we also have
\begin{equation}\label{sumsq}
x^2+y^2=2x^2+o(x^2).
\end{equation}

Using 
Taylor's decomposition,
\begin{equation}\label{ukiduki}
(1+z)^{\frac{d+1}{2}}=1+\frac{d+1}{2}z+\frac{d^2-1}{8}z^2+o(z^2),\quad |z|<1,
\end{equation}
and applying it for $q(x)$ and $q(y)$, the sum of these and (\ref{per}) results in
$$
0=\frac{d+1}{2}(q(x)+q(y))+\frac{d^2-1}{8}(q^2(x)+q^2(y))+o(q^2(x))+o(q^2(y)).
$$
This 
and (\ref{ep})  yield
\begin{multline*}
0=\frac{d+1}{2}(\varepsilon_1 x+\varepsilon_2 x^2-\varepsilon_1 |y|+\varepsilon_2 y^2)\,+\\
+\,\frac{d^2-1}{8}((\varepsilon_1 x+\varepsilon_2 x^2)^2+(-\varepsilon_1 |y|+\varepsilon_2 y^2)^2)+o(x^2)
\end{multline*}
\begin{multline*}
=\frac{d+1}{2}(\varepsilon_1 (x-|y|)+\varepsilon_2 (x^2+y^2))+\frac{d^2-1}{8}\varepsilon_1^2(x^2+y^2)+o(x^2)\\
=\frac{d+1}{2}\varepsilon_1(x-|y|)+\Big(\frac{d^2-1}{8}\varepsilon_1^2+\frac{d+1}{2}\varepsilon_2\Big)(x^2+y^2)+o(x^2).
\end{multline*}
It remains to apply   (\ref{sum}) and (\ref{sumsq}) to obtain
$$
0=\frac{d+1}{2}\varepsilon_1^2  x^2+
\Big(\frac{d^2-1}{8}\varepsilon_1^2+\frac{d+1}{2}\varepsilon_2\Big)2x^2
+o(x^2).
$$
Therefore, 
$$
\Big(\frac{d+1}{2}+\frac{d^2-1}{4}\Big)\varepsilon_1^2+(d+1)\varepsilon_2=0,
$$
or
$$
(d+1)\varepsilon_1^2+4\varepsilon_2=0.
$$
This gives the desired result by (\ref{per11}).
\ep

\begin{figure}[ht]
	\centering
	\includegraphics[height=3in]{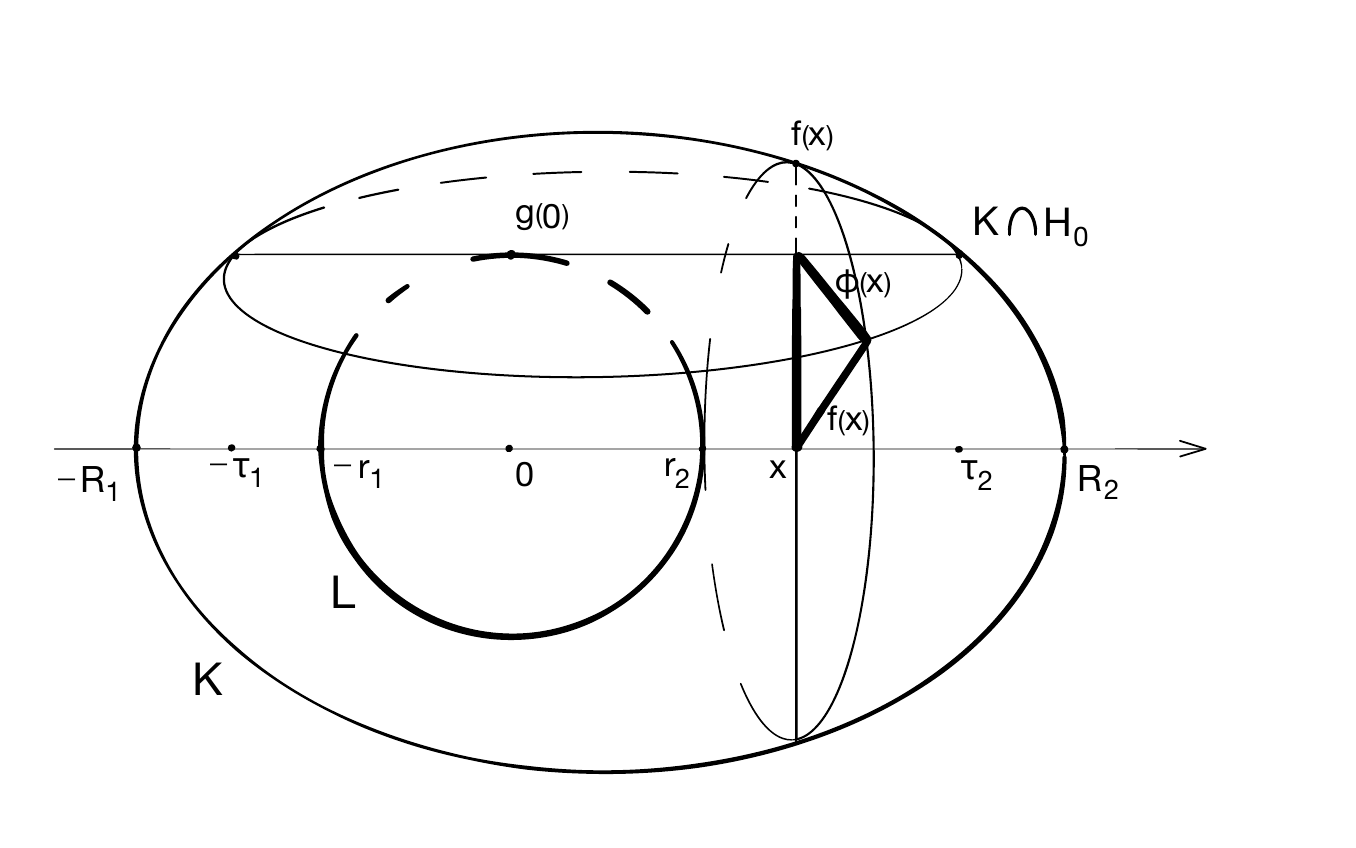} 
	\caption{We have $f^2(x)=\phi^2(x)+g^2(0)=\phi^2(x)+f^2(0)-\sigma^2$.}
	\label{Hsu}
\end{figure}

\subsection{Auxiliary formulas describing  the boundary of the horizontal section $P=K\cap H_0$} We use the notation from  the previous subsection.
If  $f$ describes the boundary of  $K$, and $\phi$ describes the boundary of the horizontal  section $K\cap H_0-(0,0,g(0))$, then
\begin{equation}\label{pip}
f^2(x)=\phi^2(x)+g^2(0)=\phi^2(x)+f^2(0)-\sigma^2,\qquad\forall x\in [-\tau_1,\tau_2],
\end{equation}
(see Figure \ref{Hsu}).

Observe  that if $\sigma$ is sufficiently close to zero, then $[-\tau_1,\tau_2]\subset [-r_1,r_2]$. On the other hand,  if $\sigma$ is large enough, then $ [-r_1,r_2]\subset[-\tau_1,\tau_2]$. The next lemma shows that in general we have only these possibilities.

\bl\label{add}
We have $ [-r_1,r_2]\subseteq[-\tau_1,\tau_2]$ or $[-\tau_1,\tau_2]\subseteq [-r_1,r_2]$. 
\el
\bp
Assume the contrary, we have
\begin{equation}\label{mozg1}
-r_1< -\tau_1,\quad r_2< \tau_2,\qquad 
\textrm{or}\quad
-\tau_1< -r_1, \quad \tau_2< r_2.
\end{equation}
We will show that  the first case  in (\ref{mozg1}) is not possible, the proof that  the second one is not possible either is similar.

 To this end, consider   the horizontal chord  inscribed into $\partial K$ and tangent to the graph of $g$ at $(0,0, g(0) )$. We have $f(-\tau_1)=f(\tau_2)=g(0)$, and $f(r_2)>f(\tau_2)$, otherwise, the points $(R_2,0,0)$, $(\tau_2, 0, f(\tau_2))$ and $(r_2,0,f(r_2))$ are on the boundary of $K$, which contradicts its convexity. On the other hand,  by (\ref{tilde}) we have $f(-r_1)=f(r_2)=\sigma$. Hence,  $f(-r_1)>f(-\tau_1)$. This contradicts the convexity of $K$, for, the points $(-R_1,0,0)$, $(-r_1,0,f(-r_1))$, and $(-\tau_1, 0, f(-\tau_1))$   must lie  on its boundary.
\ep

Let 
$$
A=s\sqrt{f^2(a(0))-a^2(s)+\chi(a(s))-\sigma^2},
$$
where for every $s\in {\mathbb R}$ we have $a(s)\in (-r_1,r_2)$.

\bl\label{compu}
Let  $s\in{\mathbb R}$, $a(s)\in (-r_1,r_2)\cap (-\tau_1,\tau_2)$ be fixed, and let 
$x$, $y$ be so small that $a(s)+x$, $a(s)-y\in (-r_1,r_2)\cap(-\tau_1,\tau_2)$. Then (\ref{chu}) and (\ref{per}) hold for these $x$, $y$, with
$$
\chi_a(x)=-2(a+A)x+\chi(a+x)-\chi(a)
$$ 
 instead of $\chi$,
and
 $q_a(x)=\frac{\chi_a(x)}{\sigma^2}$ instead of $q$.
\el

\begin{figure}[ht]
	\centering
	\includegraphics[height=3in]{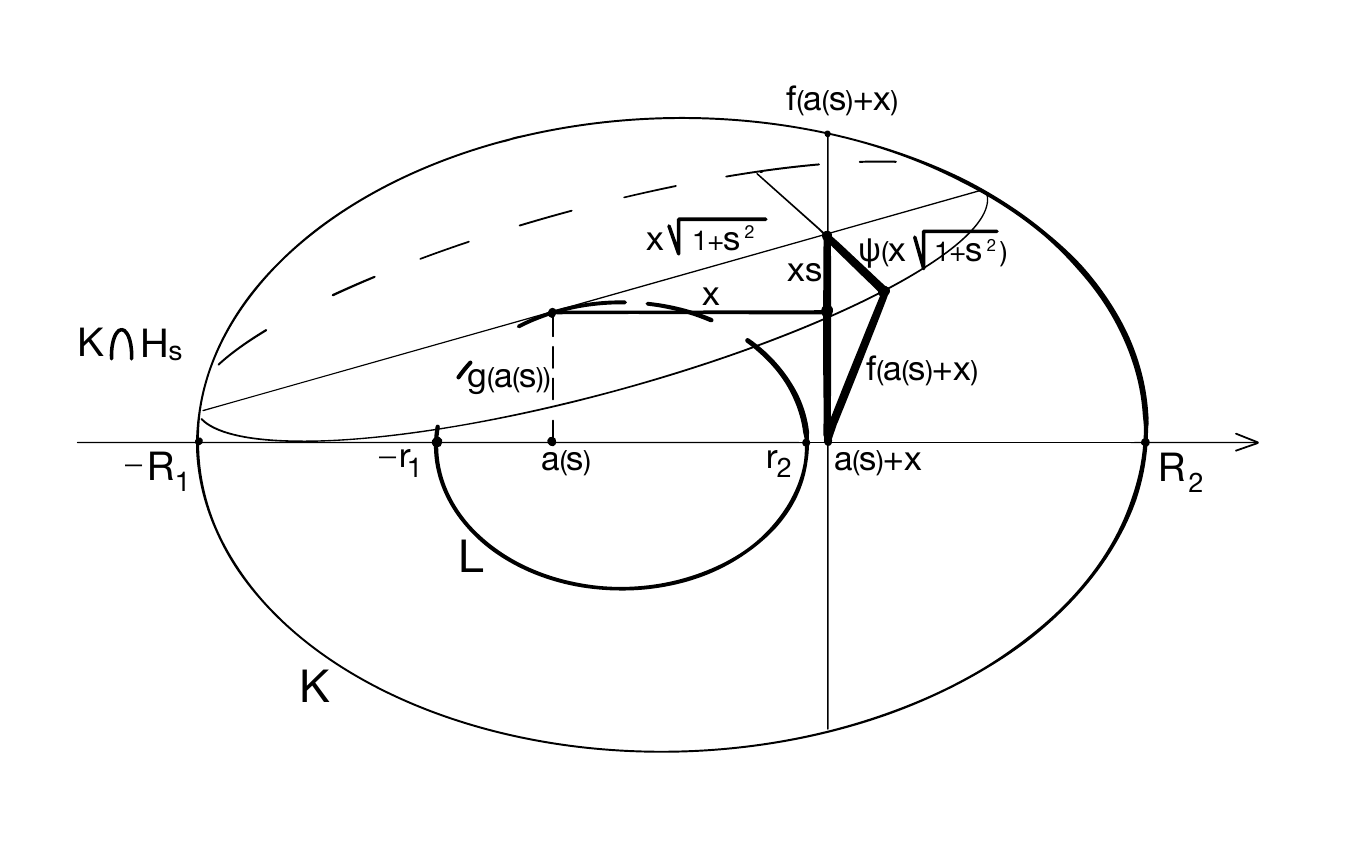} 
	\caption{We have $\psi^2(x\sqrt{1+s^2})=f^2(a(s)+x)-(g(a(s))+xs)^2$.}
	\label{Nsu}
\end{figure}

\bp
Fix any $s\in{\mathbb R}$  and  $a(s)\in (-r_1,r_2)\cap (-\tau_1,\tau_2)$.
We can assume that  the boundary of  $K\cap H_s$ is described by a positive function $\psi$ satisfying the $(d+1)$-equichordal property (we pick $(a(s),0,g(a(s))$ as the origin in $H_s$),
\begin{equation}\label{dusha}
(x^2(1+s^2)+\psi^2(x\sqrt{1+s^2}))^{\frac{d+1}{2}}\,+
\end{equation}
$$
+\,(y^2(1+s^2)+\psi^2(y\sqrt{1+s^2}))^{\frac{d+1}{2}}=2\sigma^{d+1},
$$
where
$$
\frac{\psi(x\sqrt{1+s^2})}{x}=\frac{\psi(y\sqrt{1+s^2})}{|y|}.
$$
This gives
\begin{equation}\label{yoxoxo}
|y|^{d+1}=\frac{(2\sigma^{d+1}-((1+s^2)x^2+\psi^2(x\sqrt{1+s^2}))^{\frac{d+1}{2}})x^{d+1}}{((1+s^2)x^2+\psi^2(x\sqrt{1+s^2}))^{\frac{d+1}{2}}},
\end{equation}
and $x$ and $y$ are so small that the conditions of the lemma are satisfied.

By the Pythagorean theorem (see Figure \ref{Nsu}), the assumption $a(0)=0$, and (\ref{pip}), we have
\begin{multline*}
\psi^2(x\sqrt{1+s^2})=f^2(a(s)+x)-(g(a(s))+xs)^2\\
=f^2(0)-\sigma^2+\phi^2(a(s)+x)-(\sqrt{f^2(a(s))-\sigma^2}+sx)^2\\
=f^2(0)-\sigma^2+\phi^2(a(s)+x)-(\sqrt{f^2(0)+\phi^2(a(s))-2\sigma^2}+sx)^2.
\end{multline*}
Therefore, using (\ref{chi}) we have
\begin{multline*}
x^2(1+s^2)+\psi^2(x\sqrt{1+s^2})\\
=x^2+f^2(0)-\sigma^2+\phi^2(a(s)+x)-(f^2(0)+\phi^2(a(s))-2\sigma^2)\,-\\
-\,2xs\sqrt{f^2(a(s))-\sigma^2}
\end{multline*}
\begin{multline*}
=x^2+\phi^2(a(s)+x)+\sigma^2-\phi^2(a(s))-2xs\sqrt{f^2(a(s))-\sigma^2}\\
=x^2+\phi^2(a(s)+x)+\sigma^2-\phi^2(a(s))-2xs\sqrt{f^2(a(0))+\phi^2(a(s))-2\sigma^2}\\
=\sigma^2-2a(s)x+\chi(a(s)+x)-\chi(a(s))-2xA.
\end{multline*}
Substituting the last expression into
 (\ref{dusha}) with    $y$  defined by  (\ref{yoxoxo}), 
we have
\begin{multline*}
(\sigma^2-2a(s)x+\chi(a(s)+x)-\chi(a(s)) -2xA    )^{\frac{d+1}{2}}\,+\\
+\,(\sigma^2+2a(s)y+\chi(a(s)-y)-\chi(a(s)) +2yA   )^{\frac{d+1}{2}}=2\sigma^{d+1},
\end{multline*}
where
$$
|y|^{d+1}=x^{d+1}\,\frac{2\sigma^{d+1}-(\sigma^2-2(a+A)x+\chi(a+x)-\chi(a) )^{\frac{d+1}{2}}}{(\sigma^2-2(a+A)x+\chi(a+x)-\chi(a) )^{\frac{d+1}{2}}}.
$$
This gives the desired result.
\ep

\bc\label{Adif}
Let  $s\in{\mathbb R}$ be fixed and such that $a=a(s)\in  (-r_1,r_2)\cap (-\tau_1,\tau_2)$. Then 
\begin{equation}\label{concb1}
2\sigma^2\chi''(a)+(d+1)(\chi'(a)-2(a+A))^2=0,
\end{equation}
where $A$ is as in the previous lemma.
\ec
\bp
By the previous lemma, we have (\ref{per}) and (\ref{root1}) 
with $q_a$ instead of $q$ and $\chi_a$ instead of $\chi$, $\chi_a(0)=0$. This gives
(\ref{concb}) with $\chi_a$ instead of $\chi$,
which is  the desired result. 
\ep

\subsection{Consequences  of the concavity of $\chi$ on $(-r_1,r_2)\cap (-\tau_1,\tau_2)$} Our next goal 
is to show that
\begin{equation}\label{clai}
\chi(a)= 0\qquad\forall a\in[-r_1,r_2]\cap [-\tau_1,\tau_2]. 
\end{equation}
The proof of (\ref{clai}) is contained  in the following three statements.

\bl\label{heart}
Let $\lambda_1>0$, $\lambda_2>0$ be  such that $[-\lambda_1,\lambda_2]\subseteq [-\tau_1,\tau_2]$
and
\begin{equation}\label{ab}
(\sigma^2+\chi(-\lambda_1))^{\frac{d+1}{2}}+(\sigma^2+\chi(\lambda_2))^{\frac{d+1}{2}}=2\sigma^{d+1}.
\end{equation}
If $\chi\le 0$ on $[-\lambda_1,\lambda_2]$, then $\chi= 0$ on $[-\lambda_1,\lambda_2]$. In particular, 
if $\chi\le 0$ on $[-\tau_1,\tau_2]$, then 
$\tau_1=\tau_2=\sigma$.
\el
\bp
By (\ref{chi}), we have
$$
0\le -\chi(x)\le \sigma^2-x^2\le \sigma^2,\qquad x\in [-\lambda_1,\lambda_2].
$$
By (\ref{ab}) we can assume that for all $x\in (0,\lambda_2]$ and for the corresponding $y\in [-\lambda_1,0)$ we have the equality in  (\ref{chu}). 
If $\chi(x)<0$, then 
the left-hand side of this equality   is strictly less than $2\sigma^{d+1}$. Hence, $\chi=0$ on $[-\lambda_1,\lambda_2]$.

Assume now that $[-\lambda_1,\lambda_2]= [-\tau_1,\tau_2]$. Since
$$
\tau_1^{\frac{d+1}{2}}+\tau_2^{\frac{d+1}{2}}=2\sigma^{\frac{d+1}{2}},
$$
by (\ref{chi}) we have two possibilities
\begin{equation}\label{case1}
\tau_1\le \sigma, \quad\tau_2\ge\sigma,\quad\chi(\sigma)=\phi^2(\sigma)\ge 0,
\end{equation}
or 
\begin{equation}\label{case2}
\tau_1\ge \sigma,\quad\tau_2\le\sigma,\quad \chi(-\sigma)=\phi^2(-\sigma)\ge 0. 
\end{equation}

We will consider  case (\ref{case1}), the proof for (\ref{case2}) is similar. By (\ref{case1}), $\chi(\sigma)=0$. Hence, $\phi(\sigma)=0$ and  $[0,\tau_2]=[0,\sigma]$, i.e., $\tau_2=\sigma$. This gives $\tau_1=\sigma$ and the lemma is proved.
\ep

\bl\label{osh1}
Let $\chi$ satisfy (\ref{concb1}), $\chi(0)=0$,  and let  $\chi'(0)=0$. Then $\chi=0$ on $[-r_1,r_2]\cap [-\tau_1,\tau_2]$.
\el
\bp
Using 
(\ref{concb1}) we have $\chi''(a)\le 0$ for all $a\in(-r_1,r_2)\cap (-\tau_1,\tau_2)$, i.e., $\chi$ is concave down on $(-r_1,r_2)\cap (-\tau_1,\tau_2)$.
Then  using the conditions of the lemma we get $\chi\le 0$ on $[-r_1,r_2]\cap [-\tau_1,\tau_2]$. 

Now we apply Lemma \ref{heart}.

 If  $ [-\tau_1,\tau_2]\subset [-r_1,r_2]$ we put $\lambda_1=\tau_1$, $\lambda_2=\tau_2$.

Let 
$[-r_1,r_2]\subset [-\tau_1,\tau_2]$. Consider 
the maximal segment $[-\lambda_1,\lambda_2]\subseteq [-r_1,r_2]$
for which (\ref{ab}) holds. We can assume  that $\lambda_1=r_1$ and $\lambda_2\le r_2$ (the proof in the case $\lambda_2=r_2$, $-\lambda_1\ge -r_1$ is similar). By Lemma \ref{heart} we have $\chi=0$ on $[-\lambda_1,\lambda_2]$. Therefore, 
 using (\ref{tilde}),  (\ref{chi})  and (\ref{pip}) we have
 $$
  g^2(x)=f^2(x)-\sigma^2,\quad \phi^2(x)=\sigma^2-x^2,\quad f^2(x)=\sigma^2-x^2+g^2(0),
 $$
 for all $x\in [-\lambda_1,\lambda_2]$.
 We recall that $a(0)=0$. Since for all $y\in [-\lambda_1,0]$ we have
 $$
 |(0,0,g(0))-(y,  \phi(y),g(0))|=\sigma,
 $$
 by the $(d+1)$-equichordal property we also have 
 $$
 |(0,0,g(0))-(x,  \phi(x),g(0))|=\sigma\qquad \forall x\in [0, \lambda_2].
 $$
This gives  $\lambda_2=-\lambda_1$ and $L$ must be a Euclidean ball, i.e., we can assume that $r_2=\lambda_2$.
\ep

\bl\label{osh2}
We have $\chi'(0)= 0$. 
\el
\bp
We recall that  $\chi(0)=0$.
Assume the contrary, that $\chi'(0)\neq 0$. 
Let  $\chi'(0)<0$ (the proof for  the case $\chi'(0)>0$ is similar). 

By (\ref{concb1}) we can assume that $\chi$ is concave down. Hence, 
 there exists $\varepsilon>0$ such  $\chi>0$ on $ (-\varepsilon,0)$ (we recall that $a(0)=0$).

Let   $Q=\{(x,y,z):|y|\le \sigma,\, x\le 0\}$, and we recall that  $G_s=K\cap \ell_s$ is a chord  centered at $(a(s), 0, g(a(s))$, parallel to the $y$-axis,  and  inscribed into $\partial K$ 
($G_s$ is of length $2\sigma$).

By symmetry with respect to the $xz$-plane the ends of $G_s$, $s\ge 0$,  must belong to $\partial Q\cap \partial K$. 
We will show that  for
some  small $s>0$ this is not true
and, by this, will obtain a contradiction.

To this end, let $0<\varepsilon_1<\varepsilon$  be so small that
for  $\xi\in(-\varepsilon_1,0)$ we have
\begin{equation}\label{est1}
\phi^2(\xi)=\sigma^2+\chi'(0)\xi+o(\xi)>\sigma^2\quad\forall\xi \in (-\varepsilon_1,0),
\end{equation}
where $o(\xi)$ is the remainder from the Taylor decomposition of $\phi$.
This shows  that the points on the curve $\gamma_-=\{(\xi,-\phi(\xi),g(0))$:  $\xi\in (-\varepsilon_1,0)\}$
$\subset \partial K$ do not belong to $Q$.  
By the symmetry of $K\cap H_0$  with respect to the line $K\cap H_0\cap \{(x,y,z)\in{\mathbb R^3}:\,y=0\}$, the points on the curve $\gamma_+=\{(\xi,\phi(\xi),g(0)):\,\xi\in (-\varepsilon_1,0)\}\subset \partial K$
do not belong to $Q$ either. 

Define the plane set 
$${\mathcal B}=\textrm{convhull}(\gamma_-,\gamma_+)
$$
$$
=\{(\xi,y,g(a(0)))\in{\mathbb R^3}:\, -\varepsilon_1<\xi<0\,-\phi(\xi)\le y\le \phi(\xi)  \}\subset K\cap H_0,
$$
 and let 
${\mathcal A}=\textrm{convex hull}(K\cap Q, {\mathcal B})$. By convexity of $K$ we have  ${\mathcal A}\subset K$. 
We claim  that for some $s>0$ small enough,  the 
ends of $G_s$ are  not on $\partial Q\cap \partial K$, which is a contradiction.

\begin{figure}[ht]
	\centering
	\includegraphics[height=3in]{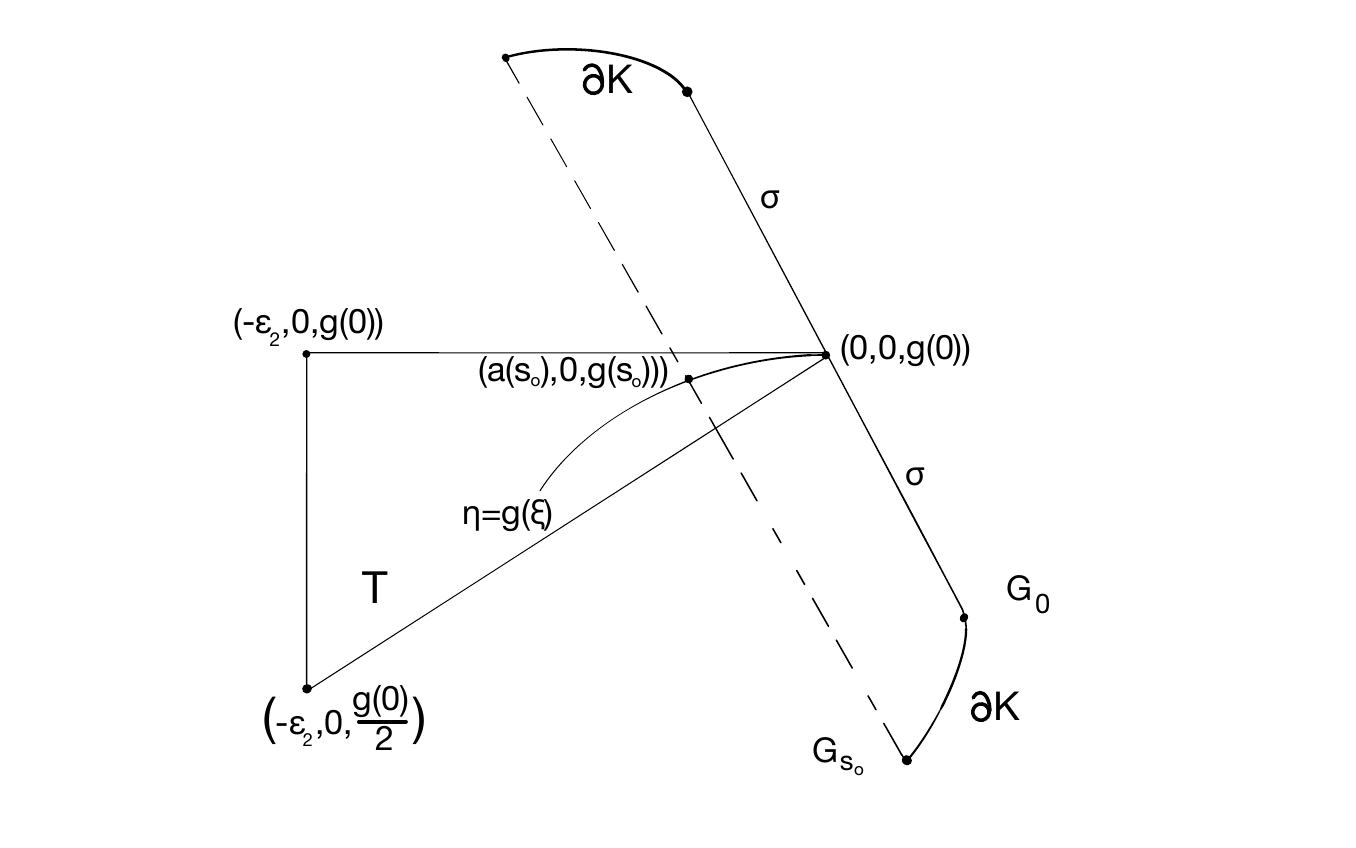} 
	\caption{The chord $G_{s_o}$ intersects $\textrm{int}\,T$, but it is longer than $2\sigma$.}
	\label{Tu}
\end{figure}

Indeed,  let $0<\varepsilon_2<\varepsilon_1$. Consider  a triangle $T$ with 
 with vertices $(0, 0, g(0))$, $(-\varepsilon_2, 0, g(0))$, $(-\varepsilon_2, 0, -\frac{g(0)}{2})$, and let $G_{(x,0,z)}$ be a chord inscribed into $\partial K$, parallel to the $y$-axis and  passing through $(x,0,z)\in \textrm{int}\,T$.
  If $\varepsilon_2$ is small enough and  $(x,0,z)\in\textrm{int}\, T$, then the ends of $G_{(x,0,z)}$
are  not in $Q$ and by (\ref{est1}) its length exceeds $2\sigma$.
Since $g$ is concave, we can find  $s_o=s_o(\varepsilon_2)>0$  so small that $a(s_o)\in (-\frac{\varepsilon_2}{2},0)$, 
and  $(a(s_o), 0, g(a(s_o)))\in \textrm{int}\,T$. In other words, the chord  $G_{s_o}$ intersects $\textrm{int}\,T$ (see Figure \ref{Tu}).
But as we noticed above,  the length of $G_{s_o}$ exceeds $2\sigma$, a contradiction. 
The proof in the case $\chi'(0)<0$ is complete.

The case $\chi'(0)>0$ can be proved similarly, one has only to consider  $a(s_o)>0$ for which $\chi(a(s_o))>0$  and  $s<0$,  and to take $Q=\{(x,y,z):|y|\le \sigma, x\ge 0\}$.
The lemma is proved.
\ep

\subsection{Conclusions}

Let $K$ and $L$ be two bodies of revolution about the $x_1$-axis in ${\mathbb R^3}$ satisfying the conditions of Theorem 
\ref{mdauzh}. We recall that  $a(0)=0$, and by  (\ref{clai})  we know  that
$\chi=0$ on $[-r_1,r_2]\cap [-\tau_1,\tau_2]$. This means that $\phi^2(x)=\sigma^2-x^2$ for all $x\in 
[-r_1,r_2]\cap [-\tau_1,\tau_2]$,  and (\ref{pip}) yields
\begin{equation}\label{ff}
f^2(x)=\phi^2(x)+f^2(0)-\sigma^2=f^2(0)-x^2
\end{equation}
for all $x\in [-r_1,r_2]\cap [-\tau_1,\tau_2]$.
Moreover, (\ref{tilde}) and (\ref{ff}) yield
\begin{equation}\label{syuzhet}
g(x)=\sqrt{f^2(0)-x^2-\sigma^2}
\end{equation}
for all $x\in [-r_1,r_2]\cap [-\tau_1,\tau_2]$, and
$$
g(-r_1)=g(r_2)=0,\quad r_1=r_2=\sqrt{f^2(0)-\sigma^2},
$$ 
provided
$[-r_1,r_2]\subseteq [-\tau_1,\tau_2]$.

\section{Auxiliary statements, the    versions of Theorem 1 of Barker and Larman,  \cite[pgs. 83-84]{BaL}}

\bl\label{iter1}
Let  $K\subset {\mathbb R^3}$ and $L\subset{\mathbb R^3}$ be two convex bodies of revolution about the $x$-axis. Assume   as above that their boundaries are described by $f$ and $g$  and satisfy   (\ref{ff}) and (\ref{syuzhet}). Then $K$ and $L$ are concentric Euclidean balls  of radii $f(0)$ and $g(0)$.
\el
\bp
 Let $\Pi$ be the $xz$-plane, and let $K\cap \Pi$ and $L\cap \Pi$ be the corresponding sections.
Observe that since $K$ and $L$ are the bodies of revolution, the sets $K\cap\Pi$ and $L\cap \Pi$ are symmetric with respect to the $x$-axis.

We will set up a certain $2$-dimensional {\it sweeping procedure} in which  the ends of the chords, that are tangent to the circular part of $\partial L\cap \Pi$ and inscribed into $\partial K\cap\Pi$, will sweep out the corresponding circular arcs on $\partial K\cap\Pi$.
Then, we will show that these   arcs  comprise 
$\partial K\cap\Pi$, thus concluding that $K\cap \Pi$ and $L\cap \Pi$  are concentric discs.

\begin{figure}[ht]
	\centering
	\includegraphics[height=3in]{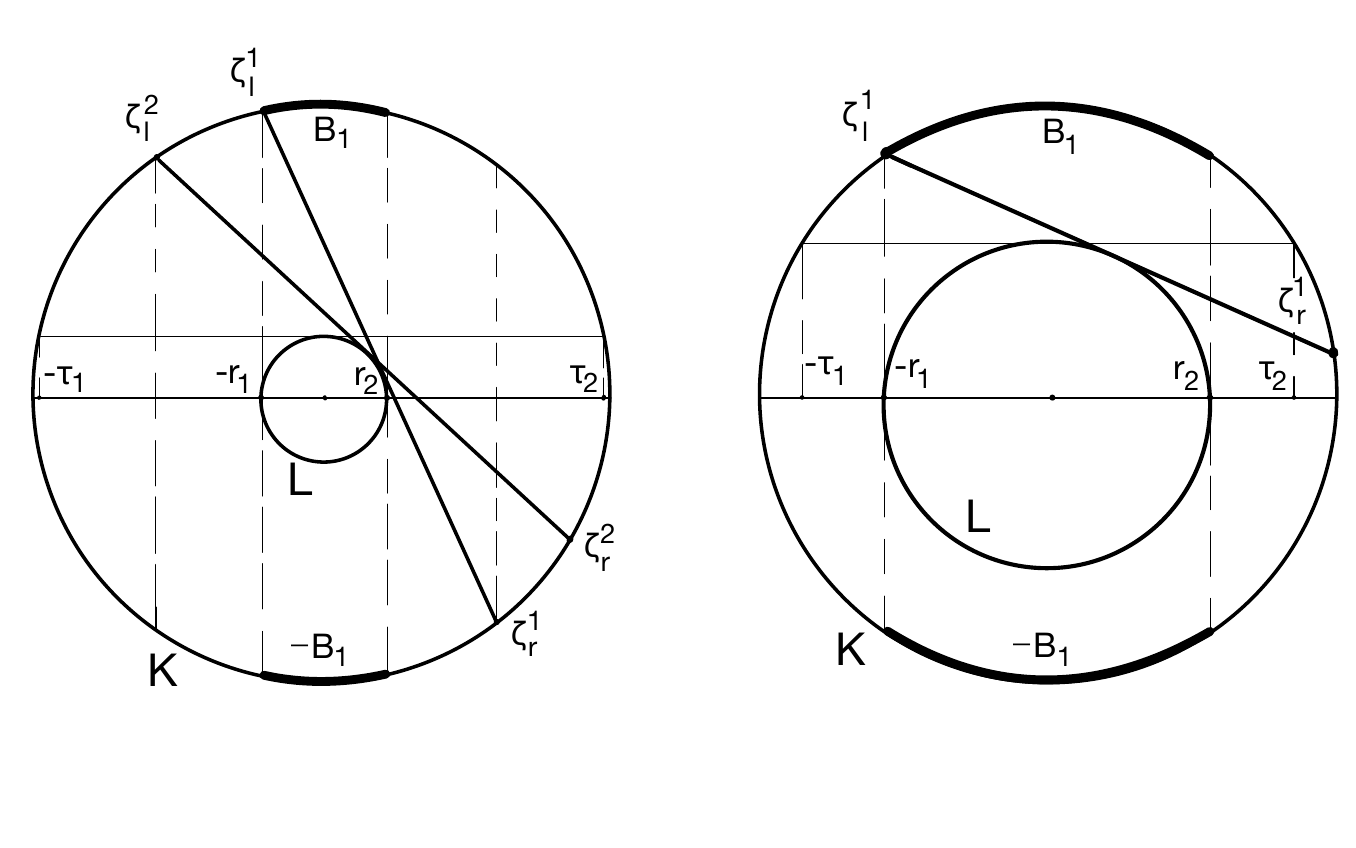} 
	\caption{The first steps in the sweeping procedure in the case  $[-r_1,r_2]\subseteq[-\tau_1,\tau_2]$. On the left $\zeta_r^1$ is below the $x$-axis, and on the right  it is above the $x$-axis. }
	\label{netv}
\end{figure}

{\it Case 1: $[-r_1,r_2]\subseteq[-\tau_1,\tau_2]$}.
As we just mentioned,    $L\cap\Pi$ is a disc of radius $g(0)=\sqrt{f^2(0)-\sigma^2}=r_1=r_2$. We will show that
$K\cap \Pi$ is a disc of radius $f(0)$.

Let ${\mathcal J}_1=[\zeta_l^1,\zeta_r^1]$ be the chord inscribed into $\partial K\cap \Pi$ and tangent to $\partial L\cap \Pi$ at $(a_1, g(a_1))$, and such that its left end is $\zeta_l^1=(-g(0), f(-g(0))$, and right end is $\zeta_1^r=(b_1, d_1)$.
We have two possibilities,
$d_1=f(b_1)>0$ or $d_1=-f(b_1)<0$ ($\zeta_r^1$ is below or  above the $x$-axis, see  Figure \ref{netv}).
Consider the arcs 
$$
{\mathcal B}_1=\{(a,f(a)):\,a\in [-g(0),g(0)]\},\quad
\beta_1=\{(a,g(a)):\,a\in [-g(0), g(0)]\},
$$
of concentric  circles, and let ${\mathcal J}^1(b)=[\zeta_l^1(b),\zeta_r^1(b)]$ be the chord
inscribed into $\partial K\cap \Pi$ and tangent to $\partial L\cap \Pi$ at $(b, g(b))\in\beta_1$ for 
$b\in [a_1, g(0)]$.
Since
the distance between $(b,g(b))$ and $\zeta_l^1(b)$ is $\sigma$, ${\mathcal J}^1(b)$ is divided by $(b,g(b))$ into two parts of equal length. 
Hence, 
while the left end  of ${\mathcal J}^1(b)$ is sweeping out  ${\mathcal B}_1$ by moving from $(-g(0), f(-g(0))$ to $(g(0),f(g(0)))$, its right end   must move  along the arc of the circle of radius $f(0)$ (from $(b_1, d_1)$ to $(g(0),-f(g(0)))$) joining 
$-{\mathcal B}_1$ from the right at $(g(0),-f(g(0)))$. 

Let $d_1=f(b_1)>0$ ($\zeta_r^1$ is  above the $x$-axis,           see the right part of Figure \ref{netv}).  Then,   the right end of 
${\mathcal J}^1(b)$ for $b\in [a_1, g(0)]$ sweeps 
out the circular part of  $\partial K\cap \Pi$ containing the one joining $(f(0),0)$ with $(g(0),-f(g(0)))$. By the aforementioned symmetry of  $\partial K\cap \Pi$  with respect to the $x$-axis, we see that the part of 
$\partial K\cap \Pi$ lying in the right half-plane is circular. Since the above procedure is symmetric with respect to the $z$-axis (we could start with the chord    ${\mathcal J}_1$ tangent to $\partial L\cap \Pi$ at $(-a_1, g(-a_1))$ and follow the sweeping arc joining  $(-b_1, d_1)$ to $(-g(0),-f(-g(0)))$), we conclude that
$\partial K\cap \Pi$
is a circle of radius $f(0)$.

Now let $d_1=-f(b_1)<0$ (see the left part of Figure \ref{netv}). By the symmetry, four points $(\pm b_1,\pm f(b_1))$ are on $\partial K\cap\Pi$ and we recall that  $f(x)=\sqrt{f^2(0)-x^2}$ for $x\in [-b_1,b_1]$.
We will repeat the above procedure for the chord 
${\mathcal J}_2$  inscribed into $\partial K\cap \Pi$ and tangent to $\partial L\cap \Pi$ at $(a_2, g(a_2))$, $0<a_2<a_1$, where 
${\mathcal J}_2=[\zeta_l^2,\zeta_r^2]$, $\zeta_l^2=(-b_1, f(-b_1))$, $\zeta_r^2=(b_2, d_2)$, 
 $b_2>b_1$, and the arcs
${\mathcal B}_2=\{(a,f(a)):\,a\in [-b_1,b_1]\}$, $
\beta_1$.

We have
two possibilities again,
$d_2=f(b_2)>0$, and $d_2=-f(b_2)<0$. 
If $d_2=f(b_2)>0$, arguing as above, we see that the part of 
$\partial K\cap \Pi$ lying in the right half-plane is circular, and by the symmetry, $\partial K\cap \Pi$ is a circle.
If $d_2=-f(b_2)<0$, taking into account that 
$(\pm b_2,\pm f(b_2))$, are on $\partial K\cap\Pi$
and $f(x)=\sqrt{f^2(0)-x^2}$ for $x\in [-b_2,b_2]$, we repeat the procedure again. 
producing the chords 
$
{\mathcal J}_3=[\zeta_l^3,\zeta_r^3]$, $0<a_3<a_2$, and etc.

If for some $j\ge 3$ we have  $d_j=f(b_j)>0$, we finish as above. 
If, on the other hand, $d_j=-f(b_j)<0$ for $j=3,4,\dots$, we produce a sequence of segments $\{[-b_j,b_j]\}_{j=1}^{\infty}$ such that $[-b_j,b_j]\subset [-b_{j+1},b_{j+1}]$, and such that $f(x)=\sqrt{f^2(0)-x^2}$ for $x\in [-b,b]$,  $b=\lim\limits_{j\to\infty}b_j$.

We can also assume that $d_j=-f(b_j)<-g(0)$ for all $j=3,4,\dots$. Indeed,   since the points $(\pm b_j,\pm f(b_j))$ must be on  $\partial K\cap \Pi$, then the condition 
$-f(b_j)\ge -g(0)$ for some $j\ge 3$ implies that
the chord with its left end at 
$(-b_j, f(-b_j))$ must have a positive second coordinate  for its right end, so
 $d_{j+1}=f(b_{j+1})>0$.

We claim that  $\partial K\cap \Pi$ is a circle.
Indeed, let $b<f(0)$ (otherwise, we are done). 
If  $-f(b)\ge -g(0)$, then 
the points $(\pm b,\pm f(b))$ must be on  $\partial K\cap \Pi$. Hence, the chord with its left end at 
$(-b, f(-b))$ must have $(b, f(b))$ for its right end,  $f(b)>0$, and we are done.

Finally, let  $-f(b)<-g(0)$ and let 
$$
{\mathfrak b}=\sup\{x\in [0,f(0)]:\,  f(x)=\sqrt{f^2(0)-x^2}\,\,\,\textrm{on} \,\,\,[0,{\mathfrak b}] \}.
$$
Then $-f({\mathfrak b})\ge -g(0)$,  otherwise
$(\pm {\mathfrak b},\pm f({\mathfrak b}))$ are on 
$\partial K\cap \Pi$, and we can repeat the procedure, contradicting the definition of ${\mathfrak b}$.

This finishes the proof of {\it Case 1}.

{\it Case 2: $[-\tau_1,\tau_2]\subsetneq [-r_1,r_2]$}.

\begin{figure}[ht]
	\centering
	\includegraphics[height=3in]{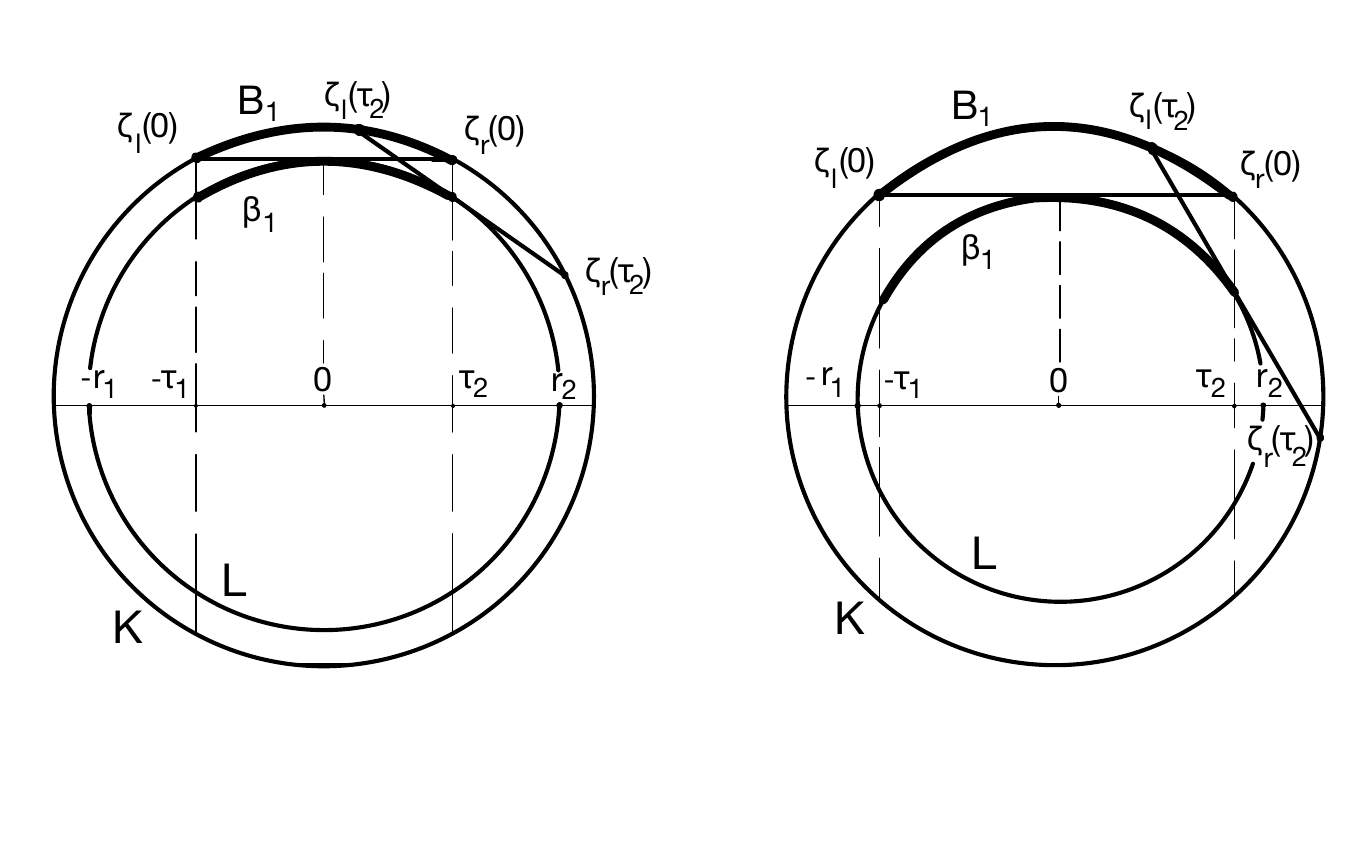} 
	\caption{The first steps in the sweeping procedure in the case  $[-\tau_1,\tau_2]\subseteq [-r_1,r_2]$. On the left $\zeta_r(\tau_2)$ is above the $x$-axis, and on the right  it is below the $x$-axis. }
	\label{CHu}
\end{figure}

Let ${\mathcal J}_a=[\zeta_l(a),\zeta_r(a)]$ be a chord inscribed into $\partial K\cap \Pi$ and tangent to $\partial L\cap \Pi$ at the point $(a,g(a))$, for some $a\in (-\tau_1,\tau_2)$ (see Figure \ref{CHu}).

Consider the arcs of concentric  circles  ${\mathcal B}_1=\{(a,f(a))_{a\in [-\tau_1,\tau_2]}\}$ and $\beta_1=\{(a,g(a))_{a\in [-\tau_1,\tau_2]}\}$ and observe that for any $b\in [0, \tau_2]$ the distance between $(b,g(b))$ and $\zeta_l(b)$ is $\sigma$ and ${\mathcal J}_b$ is divided by $(b,g(b))$ into two parts of equal length. 
Hence, 
while the left end $\zeta_l(b)$ for $b\in [0, \tau_2]$  is sweeping out the  part of ${\mathcal B}_1$ by moving from $\zeta_l(\tau_2)$ to $\zeta_l(0)=(-\tau_1,f(\tau_1))$, the right end $\zeta_r(b)$  must move  along the arc of a circle of radius $f(0)$ joining 
${\mathcal B}_1$ from the right at $\zeta_r(0)$. 
If we denote the coordinates of $\zeta_r(\tau_2)$ by $(v_1, d_1)$, we 
have
(\ref{ff}) and (\ref{syuzhet})  for all $x$ in the interval 
 $[-\tau_1, v_1]$ strictly containing $[-\tau_1, \tau_2]$. 

As in  {\it Case 1} we have two possibilities, $d_1=f(v_1)<0$ or $d_1=-f(v_1)<0$.

If $d_1=-f(v_1)<0$, we stop the procedure and see that the parts of $\partial K\cap \Pi$ and $\partial L\cap \Pi$, located in the right half-plane are concentric circles.

Let $d_1=f(v_1)>0$. Denote by $\alpha\in (0,\frac{\pi}{2})$ the angle
between the tangent line to $\partial L\cap\Pi$ passing through $(R_2,0)$ and the $x_1$-axis (we recall that $(R_2,0)$ is the point of intersection of $\partial K$ with the $x$-axis).  If
$\alpha_1$ is the angle between the line containing
${\mathcal J}_{\tau_2}$ and the $x_1$-axis,
then  
 $v_1-\tau_2=\sigma\cos\alpha_1$, 
and by convexity $\alpha_1<\alpha$. 
We repeat the process with  the larger arcs of concentric  circles  ${\mathcal B}_2=\{(a,f(a))_{a\in [-\tau_1,v_1]}\}$ and  $\beta_2=\{(a,g(a))_{a\in [-\tau_1,v_1]}\}$
instead of  ${\mathcal B}_1$ and $\beta_1$. 
As above we have two possibilities 
$d_2=f(v_2)>0$ or $d_2=-f(v_2)<0$ for  the corresponding right end $(v_2, d_2)$ of the chord ${\mathcal J}_{v_1}$. If $d_2=-f(v_2)\le 0$, we stop. If 
$d_2=f(v_2)>0$ we repeat, observing that
$v_2-v_1=\sigma\cos\alpha_2$ for the angle  $\alpha_2$ between 
the line containing
${\mathcal J}_{v_1}$ and the $x_1$-axis,
$\alpha_2<\alpha$.
Proceeding  this way, we will construct 
 the corresponding arcs ${\mathcal B}_j$ and $\beta_j$,   $j=3,\dots, m$.
If for some $j$ we have $d_j=-f(v_j)<0$ for the corresponding right end of the chord ${\mathcal J}_{v_{j-1}}$, we will stop. Otherwise, we will proceed with $d_j=f(v_j)>0$ for all $j=3,\dots, m$, and the corresponding angles $\alpha_j<\alpha$.
Since $v_j-v_{j-1}=\sigma\cos\alpha_j\ge \sigma\cos\alpha$ for $j=2,\dots,m$,
we will have
$$
v_m=\tau_2+(v_1-\tau_2)+\dots+v_m-v_{m-1}
\ge \tau_2+m\sigma\cos\alpha\ge R_2,
$$
provided $m$ is large enough. 
We have proved that the parts of $\partial K\cap \Pi$ and $\partial L\cap \Pi$, located in the right half-plane are concentric circles.

Similarly,  while the right end $\zeta_r(b)$ for $b\in [-\tau_1,0]$  is sweeping out the  part of ${\mathcal B}_1$ by moving from  $\zeta_r(-\tau_1)$ to $\zeta_r(0)=(\tau_2, f(\tau_2))$, the left end $\zeta_l(b)$  must move  along the arc of a circle of radius $f(0)$ joining 
${\mathcal B}_1$ from the left at $\zeta_l(0)$. 
If we denote the coordinates of $\zeta_l(-\tau_1)$ by $(-u_1, {\mathfrak d}_1)$, we 
have
(\ref{ff}) and (\ref{syuzhet})  for all $x$ in  the interval $[-u_1,\tau_2]$ strictly containing $[-\tau_1, \tau_2]$. 
This gives
(\ref{ff}) and (\ref{syuzhet})  for all $x\in [-u_1,\tau_2]$.
Considering two cases ${\mathfrak d}_2=f(\tau_2)>0$ or ${\mathfrak d}_2=-f(\tau_2)<0$, we can repeat the argument above to obtain that $\partial K\cap \Pi$ and $\partial L\cap \Pi$ are concentric discs.
\ep

\bl\label{vm}
Let $K$ and $L$ be two convex bodies in ${\mathbb R^3}$  satisfying the conditions of Theorem \ref{mdauzh}. If $L$ is a body of revolution, then $K$ is also a body of revolution with the same axis of rotation.
\el
\bp
We assume that the $x$-axis is the axis of rotation of $L$. 
We will set up a $3$-dimensional sweeping procedure  rotating  the cones that are  tangent to $\partial L$ with vertices on $\partial K$. 

Let $W_x$ be a plane parallel to the $yz$-plane and passing through $(x,0,0)$, $x\in {\mathbb R}$, and let 
$M(x)\subset W_x$ be a circle centered at $(x,0,0)$.
We will show that for every $x$ such that $(x,y,z)\in \textrm{int}K$, the 
generators of the sweeping  cones cut out the circles $M(x)\subset\partial K$, 
thus proving that $K$ is a body of revolution about the $x$-axis.

Let  $e'=(x',0,0)$,  $e''=(x'',0,0)$,  be two points of the  intersection of the $x$-axis with $\partial K$, $x'>0$,  $x''<0$. 
To set up the procedure, we will make several auxiliary remarks and observations.

\begin{figure}[ht]
	\centering
	\includegraphics[height=3in]{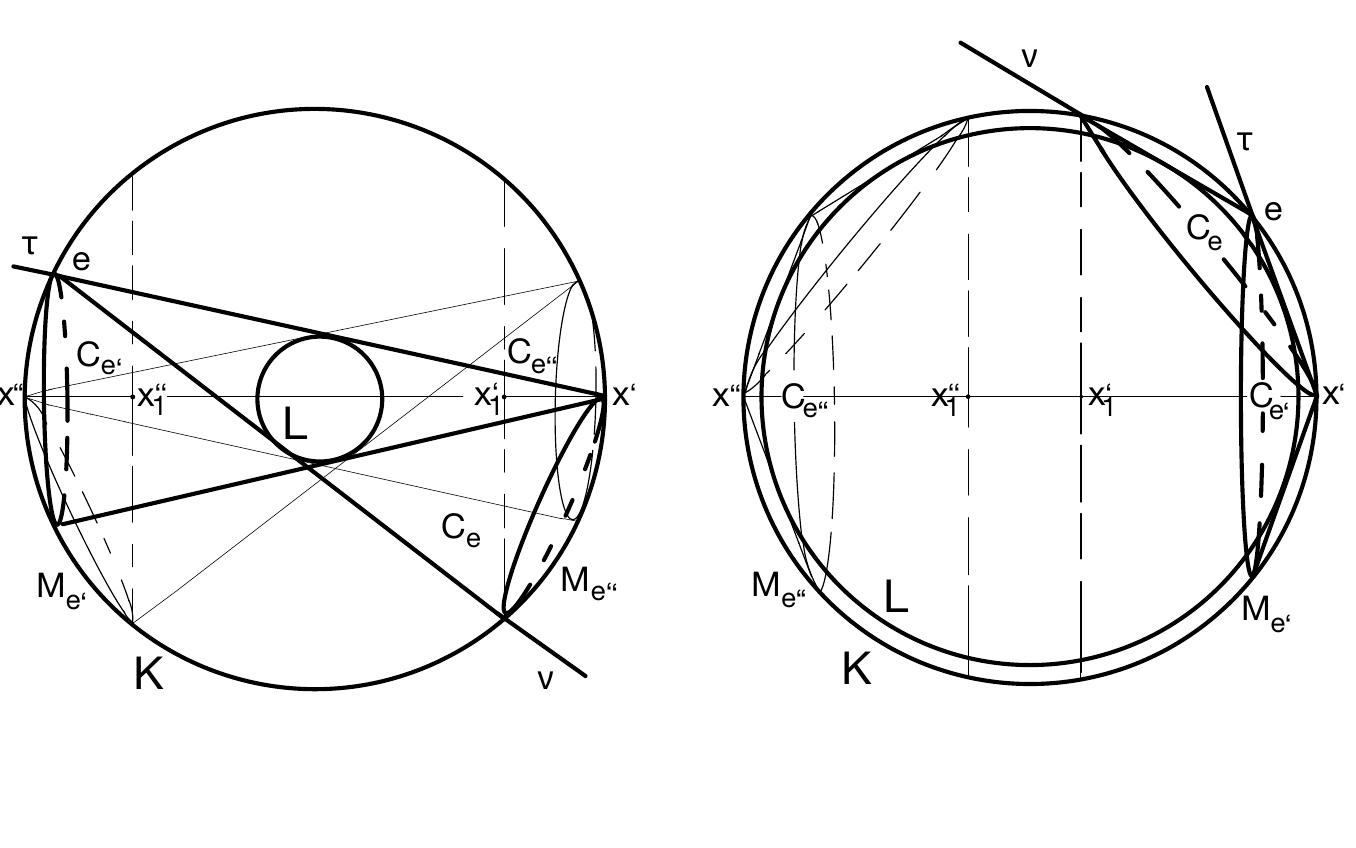} 
	\caption{The sweeping cones in the $3$-dimensional procedure. On the left part we have $M_{e'}$ is left  to $M_{e''}$ and on the right part $M_{e'}$ is  right  to $M_{e''}$.  }
	\label{rabota1}
\end{figure}

By the $(d+1)$-equatorial property of $K$ and $L$, for every ray $\tau$ emanating from $e'$ and tangent to $\partial L$ we have
\begin{equation}\label{vl}
|e'-\partial L\cap \tau|^{d+1}+|\partial L\cap \tau-\partial K\cap \tau|^{d+1}=2\sigma^{d+1}.
\end{equation}
Since $|e'-\partial L\cap \tau|$ is  constant independent of $\tau$, by (\ref{vl}) we see that the same is true 
for  $|\partial L\cap \tau-\partial K\cap \tau|$. Therefore,  for all rays $\tau$ emanating from $e'$ and tangent to $\partial L$,
all the   chords $K\cap\tau$ have the same length. 
Since $L$ is the  body of revolution, 
for any rotation  $\Phi=\Phi_{\varphi}$  by the angle $\varphi\in (0,2\pi)$ around the $x$-axis,  the points 
$\{\partial L\cap \Phi_{\varphi}\tau:\varphi\in [0,2\pi]\}$ form
a circle  centered on the $x$-axis.
By  similarity of triangles, the ends $\{\partial K\cap \Phi_{\varphi}\tau\neq e':\,\varphi\in[0,2\pi]\}$
of the chords $K\cap\Phi_{\varphi}\tau$
form a circle  $M_{e'}=\partial K\cap C_{e'}$ centered  on the $x$-axis, 
where $C_{e'}$ is the cone tangent to $\partial L$ with the vertex at $e'$ (see Figure \ref{rabota1}).

Now  we take any point $e\in M_{e'}\subset\partial K$ and repeat  a similar  argument for   the cone $C_e$ tangent to $\partial L$ with the vertex at $e$. 
Observe  that for any 
ray $\nu$ generating $C_e$, the ends $\{\partial K\cap \Phi_{\varphi}\nu\neq \Phi_{\varphi}e:\,\varphi\in[0,2\pi]\}$ of the chords $K\cap\Phi_{\varphi}\nu$ form a circle $M_{e,\nu}(e')\subset \partial K$ with the center on the $x$-axis and which is parallel to $M_{e'}$ (see Figure \ref{rabota1}).

Indeed, let $e\in M_{e'}$ and  let $\nu$ be any ray generating $C_e$. 
By rotation invariance of the length, $|\Phi(K\cap\nu)|=|K\cap\nu|$, and by the rotation invariance of  $L$, $|e-\partial L\cap\nu|=|\Phi(e)-\partial L\cap\Phi(\nu)|$.  Since
(\ref{vl}) holds with $e$, $\nu$, and $\Phi(e), \Phi(\nu)$,  instead of $e'$, $\tau$,  and since  for $\varphi\in[0,2\pi]$ the points $\Phi_{\varphi}(e)$ and $\partial L\cap\Phi_{\varphi}(\nu)$, ``move along" the circles centered on the $x$-axis, 
we see that the ends 
$\{\partial K\cap \Phi_{\varphi}\nu\neq \Phi_{\varphi}e:\,\varphi\in[0,2\pi]\}$ of the chords $K\cap\Phi_{\varphi}\nu$
form a circle $M_{e,\nu}(e')$ parallel to $M_{e'}$ and centered on the $x$-axis. This proves the observation.

We can repeat the same argument with $e''$ instead of $e'$.

Now we are ready to make the first step of our procedure. Let 
$$
x_1'=\inf\{x:\, M(x)=M_{e,\nu}(e')  \quad  \textrm{with}\,\,\,\nu\,\,\, \textrm{generating}\,\,\, C_e, \,e\in M_{e'}  \},
$$
$$
 x_1''=\sup\{x:\, M(x)=M_{e,\nu}(e'')  \quad  \textrm{with}\,\,\,\nu\,\,\, \textrm{generating}\,\,\, C_e, \,e\in M_{e''}  \}.
$$
Observe that $ x_1'<x'$ and 
$x''<x_1''$.

We will consider two cases, $M_{e'}$ is right to $M_{e''}$ and $M_{e'}$ is  left  to $M_{e''}$ (see Figure \ref{rabota1}).
In both cases, by the above observations we have   $M(x)\subset\partial K$ $\forall x\in [x_1',x']\cup [x'',x_1'']$,
 i.e.,  the sets $\{(x,y,z)\in K:$  $x\in [x_1',x']\}$ and $\{(x,y,z)\in K:$  $\cup [x'',x_1'']\}$
are the bodies of revolution about the $x$-axis.

\begin{figure}[ht]
	\centering
	\includegraphics[height=3in]{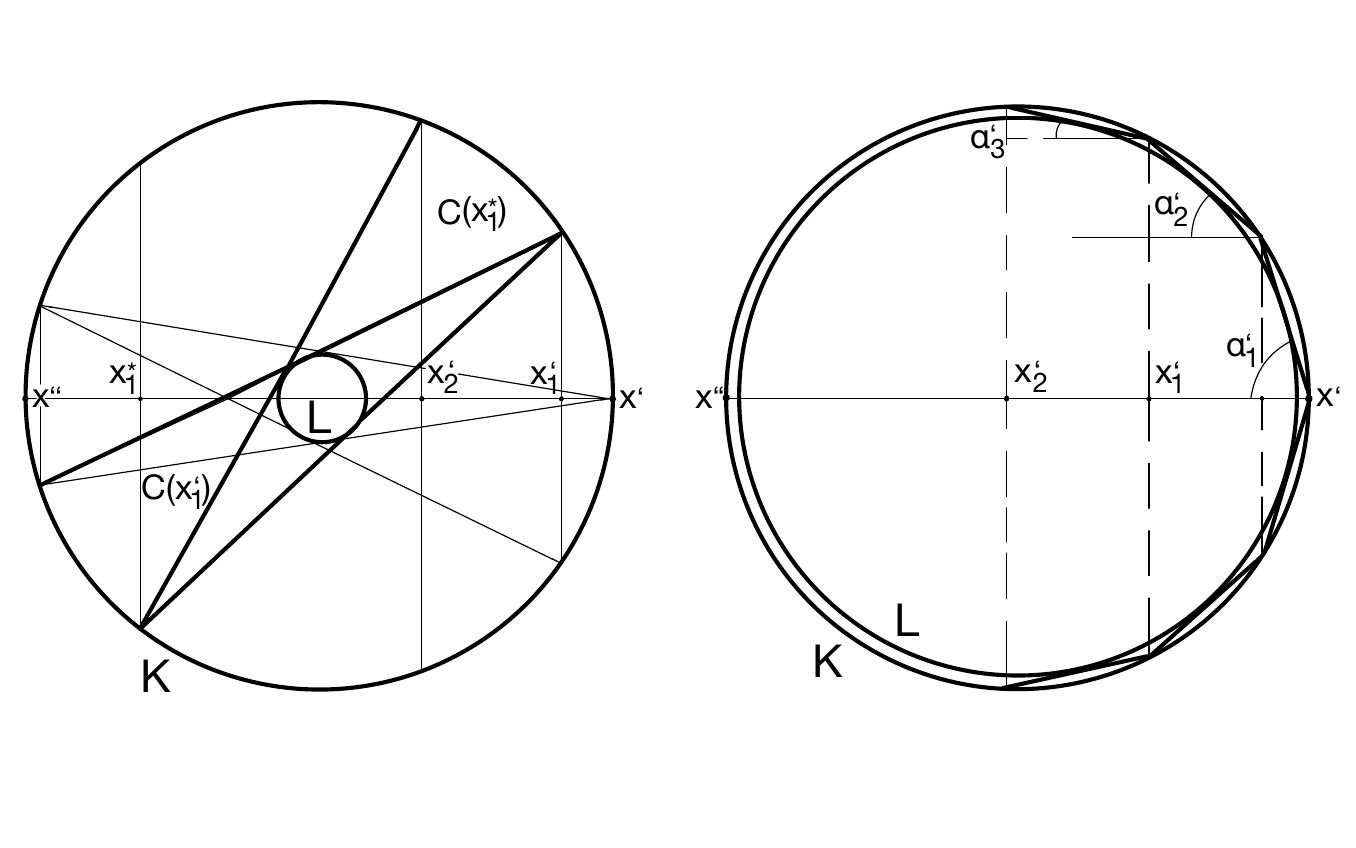} 
	\caption{Sections of the  sweeping cones by the $xz$-plane and the  construction of $x_1'$ and  $x_2'$. }
	\label{rabota2}
\end{figure}

Let $M_{e'}$ be right to $M_{e''}$.
We repeat the above argument for the generators of the cone $C_e$, with $e$ belonging to the circles $M(x_1')$ and  $M(x_1'')$ . This gives 
$M(x)\subset\partial K$ $\forall x\in
[x_2',x_1']$ for some $x_2'\in[x'',x_1')$ (see the right part of Figure \ref{rabota2}), and, similarly,  $M(x)\subset\partial K$ $\forall x\in[x_1'',x_2'']$ for some $x_2''\in (x_1'',x']$, and etc.

We claim that after $m\in {\mathbb N} $ steps  we have 
$M(x)\subset\partial K$ $\forall x\in [x'',x']$, i.e., $K$ is a body of revolution.
In fact, since
 the lengths of all chords tangent to $\partial L$ and inscribed into $\partial K$ exceed or equal to 
 $$
 2^{\frac{1}{d+1}}\sigma=\min\{x+y:\,x^{d+1}+y^{d+1}=2\sigma^{d+1}\quad \textrm{and}\quad x\ge 0, y\ge 0\},
 $$
 we have
$x_j'-x_{j+1}'\ge 2^{\frac{1}{d+1}}\sigma\cos\alpha_j'$.
Here by convexity $\alpha_{j+1}'<\alpha_{j}'<\frac{\pi}{2}$ for $j=0,1,\dots,m$, $x_0'=x'$, $x_0''=x''$
(see the right part of Figure \ref{rabota2}).
Similarly, $x_{j+1}''-x_j''\ge 2^{\frac{1}{d+1}}\sigma\cos\alpha_j''$ for the corresponding $\alpha_j''$.
Hence,  for  sufficiently large $m$ we have
$$
\sum\limits_{j=0}^m((x_j'-x_{j+1}')+(x_j''-x_{j+1}''))
\ge \sum\limits_{j=0}^m2^{\frac{1}{d+1}}\sigma(\cos\alpha_j'+\cos\alpha_j'')\ge x'-x'',
$$
and the claim is proved.

It remains to consider the case where $M_{e'}$ is left to $M_{e''}$. 
 As above, we  will run the procedure that starts at  $e'$ and follows the cones $C_e$ with vertices  $e$ at $M(x_1')$, $M(x_2')$, $\dots$, $M(x_m')$. This time, however, each point $x_j'$, $j=2$, $\dots$, $m$, will be constructed slightly differently:
for  the cones $C(x_1')=C_e(x_1')$, tangent to $\partial L$ with $e\in M(x_1')$, define
$$
x_1^*=\sup\{x:\,(x,y,z)\in (C_e(x_1')\cap \partial K)\setminus M(x_1')\quad\textrm{for}\quad e \in M(x_1')\};
$$ 
in its turn,  for  the cones $C(x_1^*)=C_e(x_1^*)$, tangent to $\partial L$ with $e\in M(x_1^*)$, let
$$
x_2'=\inf\{x:\,(x,y,z)\in (C_e(x_1^*)\cap \partial K)\setminus M(x_1^*)\quad\textrm{for}\quad e \in M(x_1^*)\},
$$
(see the left part of Figure \ref{rabota2}). 
Observe that $x_2'<x_1'$,  and for all $x\in [x_2',x_1']$ we have
 $M(x)\subset \partial K$.

We can repeat the construction with  the  corresponding $x_j^*$ and $x_j'$, $j=2,\dots,m$,   to see that 
$M(x)\subset \partial K$ for $x\in [x_m',x']$. 
Let $r_1$ and $r_2$ be such that $\{x:\,(x,y,z)\in L    \}$ $=$ $[-r_1,r_2]$.
If $0<\inf\limits_{j\ge 2}x_j'\le r_2$, we stop the procedure. For, considering the cone $C_e$ tangent to $\partial L$ with $e\in M(\inf\limits_{j\ge 2}x_j')$, we see that the parts of $K$ and $L$ in $\{(x,y,z)\in
	{\mathbb R^3}:\, x\ge 0\}$ are bodies of revolution.

Assume now that $\inf\limits_{j\ge 2}x_j'> r_2$, and let $$
\gamma=\inf\{x:\, M(y)\subset \partial K        \quad\forall \,y\ge x\},\quad\quad  0\le\gamma\le \inf\limits_{j\ge 2}x_j',
$$
(without loss of generality we can assume that $\gamma>r_2$, otherwise we finish as above).
If $\sup\limits_{j\ge 2}x_j^*\ge -r_1$, we stop. In this case, considering the cone $C_e$ tangent to $\partial L$ with $e\in M(\sup\limits_{j\ge 2}x_j^*)$, we see that the part of $K$ in $\{(x,y,z)\in
{\mathbb R^3}:\, x\ge 0\}$ is a  body of revolution. Finally, the case 
$\sup\limits_{j\ge 2}x_j^*< -r_1$ is impossible, for, we could continue the procedure, which contradicts the definition of $\gamma$.

Thus, the parts of $K$ and $L$ in $\{(x,y,z)\in
{\mathbb R^3}:\, x\ge 0\}$ are bodies of revolution.
The analogous 
argument for $\{(x,y,z)\in
{\mathbb R^3}:\, x\le 0\}$ corresponding to $e''$ follows  similarly.
\ep

\section{Proof of Theorem \ref{mdauzh}}\label{ux}

 Let $L$ be a body of revolution about the $x_1$-axis
 and let $W$ be any $3$-dimensional subspace containing the $x_1$-axis.
If $d\ge 4$, we will consider $K\cap W$, $L\cap W$,
where 
 without loss of generality we assume that $W=\{x\in{\mathbb R^d}:\,x_4=\dots=x_d=0\}$.
 
 By Lemma \ref{vm} we know that $K\cap W$ and $L\cap W$ are bodies of revolution about the $x_1$-axis.
 It follows that, by 
  Lemmas \ref{osh1} and \ref{osh2}, we have 
 (\ref{ff}) and (\ref{syuzhet}). Hence, by Lemma \ref{iter1}, $K\cap W$ and $L\cap W$ are the concentric Euclidean balls.
 
 Let now $\Pi$ be any $2$-dimensional subspace of ${\mathbb R^d}$, and let $e_1$ be the first coordinate vector.  If $e_1\notin \Pi$, let  $W_{\Pi}=\textrm{span}(\Pi, e_1)$, and if $e_1\in \Pi$ let $W_{\Pi}$ be any $3$-dimensional subspace containing $\Pi$. In both cases, by the above, $K\cap W_{\Pi}$ and $L\cap W_{\Pi}$ are the concentric Euclidean balls. Hence, $K\cap \Pi$ and $L\cap \Pi$ are the concentric discs.
 Since $\Pi$ was chosen arbitrarily, 
the application of   \cite[Corollary 7.1.4, page 272]{Ga}  finishes the proof of Theorem \ref{mdauzh}.

\vspace{1cm}

{\bf Acknowledment}. The author is very thankful to Alexander Fish and Vlad Yaskin  for their help and numerous very useful discussions. He is also very indebted to the referees who  found the gap in the original version of the paper and made many suggestions that helped to improve it.

\end{document}